\documentclass[11pt]{article}
\usepackage{amsmath,amsfonts,amssymb,bm}

\usepackage{fancyhdr}
\font\Ttl=cmb10 scaled \magstep2
\font\Srm=cmr9 
 
\font\Sbf=cmb9
\font\bfs=cmb10 scaled \magstep1

\def\RR{\mathbb{R}}

\def\CC{\mathbb{C}}

\long\def\DEL#1{{}}

\begin{document}

\font\Goth=eufm10
\def\Psig{\bm{\Psi}}
\def\PPP{\hbox{\Goth P}}
\def\FFg{\hbox{\Goth F}}

\def\norm#1{\big\Vert #1 \big\Vert}
\def\abs#1{\big\vert #1 \big\vert}

\def\eps{\varepsilon}

\def\Chi{{\bm X}}

\def\Tto{\rightrightarrows}



\def\MDPI{{1}}
\def\STA{{2}}
\def\Zlam{{3}}
\def\Serg{{4}}
\def\Hahm{{5}}
\def\Cao{{6}}
\def\WHIT{{7}}
\def\Berger{{8}}

\newcounter{SEC}
\newcounter{THM}
\def\ADS{\addtocounter{SEC}{1} \setcounter{THM}{0}}
\def\ADT{\addtocounter{THM}{1}}
\def\NUM{{\theSEC.\theTHM}}

\lhead{\textsc{L.L. Stach\'o}} 
\rhead{$C_1$-splines over triangular meshes}

\centerline{\Ttl On Hermitian interpolation of first order data with 
$\phantom{\displaystyle{\int_0^1}}$}
\centerline{\Ttl locally generated C1-splines over triangular meshes}
\bigskip
\centerline{\textsc{L.L. STACH\'O}}

\bigskip\bigskip
{\narrower \Srm \noindent
{\Sbf Abstract}. 
Given a system of triangles in the plane $\RR^2$
along with given data of function and gradient values
at the vertices, 
we describe the general pattern of 
local linear methods invoving only four smooth standard
shape functions which results in a spline function
fitting the given value and gradient data value with
${\cal C}^1$-coupling along the edges of the triangles.
We characterize their invariance properties
with relavance for the construction of interpolation surfaces
over triangularizations of scanned 3D data. 
The numerically simplest procedures among them 
leaving invarant all polynomials of 2-variables with degree 0 resp 1
involve only 
polynomials of 5-th resp. 6-th degree,
but the characteizations give rise to a huge variety of procedures 
with non-polynomial shape functions.
\par}
\renewcommand{\thefootnote}{\fnsymbol{footnote}}
\footnote[0]{\hskip-6mm
Received March 29, 2023. \hfill\break
\hfill\break
}

\vskip10mm
\ADS
\noindent{\bfs \theSEC. 
Introduction}

\bigskip
\rm
Recently [\MDPI] 
we published a ${\cal C}^1$-spline interpolation method
over 2-dimensional triangular meshes
with  polynomials of 5-th degree
with low operational costs:  
using first order data (value and gradient) at the mesh vertices
for input.
The result relies on the following
basic tool: given  a non-degenerate triangle with vertices 
${\bf p}_1,{\bf p}_2,{\bf p}_3 \in\RR^2$ along with three
values $f_1,f_2,f_2\in\RR$ respectively three
linear functionals $A_1,A_2,A_3 \in{\cal L}(\RR^2,\RR)$
and three vectors ${\bf u}_1,{\bf u}_2,{\bf u}_3 \in\RR^2$
such that ${\bf u}_k \not\parallel {\bf p}_i-{\bf p}_j$,
we can construct a polynomial $F:\RR^2\to\RR$ of 5-th degree  
of the form \ $F = F_0 - H$ \ where 
\ADT
\newcounter{FS} \setcounter{FS}{\value{SEC}}
\newcounter{FT} \setcounter{FT}{\value{THM}}
$$\begin{aligned}
&F_0({\bf x})= \sum_{i=1}^3 \Big[
\Phi\big(\lambda_i({\bf x})\big) f_i + 
\Theta\big(\lambda_i({\bf x})\big) A_i({\bf x}-{\bf p}_1) \Big] ,\cr 
&H({\bf x}) = 30\lambda_1({\bf x})^2\lambda_2({\bf x})^2\lambda_3^2({\bf x})^2
\sum_{k=1}^3 
 \lambda_k({\bf x})^{-1} {M_{k} {\bf u}_k \over G_k {\bf u}_k}
\end{aligned}
\leqno(\NUM)$$
in terms of the {\it barycentric weights} 
$\lambda_1,\lambda_2,\lambda_3:\RR^2\to\RR$
of the triangle 
${\bf T} = {\rm Conv}\{ {\bf p}_1,{\bf p}_2,{\bf p}_3\}$
associated with the vertices 
and the {\it shape functions} 
\ADT
\newcounter{SHS} \setcounter{SHS}{\value{SEC}}
\newcounter{SHT} \setcounter{SHT}{\value{THM}}
$$\Phi(t) = t^3 (10-15t+6t^2), \qquad \Theta(t)=t^3(4-3t)
\leqno(\NUM)$$
for a Hermite interpolation $F_0$ on ${\bf T}$ with the data
$f_i,A_i$ $(i=1,2,3)$.
The correction term $H$ is defined 
by means of the linear functionals 
$\RR^2\to\RR$
$$\begin{aligned}
&G_k {\bf u} = \lambda_k^\prime {\bf u} = 
{d\over dt}\Big\vert_{t=0} \lambda_k ({\bf x} + t{\bf u}) =
 \lambda_k ({\bf x} + {\bf u}) -  \lambda_k ({\bf x}), 
\cr 
&M_{k} {\bf u} = 
\sum_{
\{ i,j \}\!=\{ 1,2,3\}\!\setminus\!\{ k\}}
[ G_i {\bf u} ] \big(  30 f_i + 12 A_i({\bf p}_j - {\bf p}_i ) \big) .
\end{aligned}$$
It is worth to notice that
the 
function $F_0$ 
in $(\theFS.\theFT)$ 
restricted to the edge $[{\bf p}_i,{\bf p}_k]$
depends only on the terms $f_i,f_j,A_i,A_j$ and the impact of
the correction by adding $H$ results in the {\it reduced side derivatives} 
({\bf RSD})
$$F^\prime({\bf x}_t) {\bf u_k} = 
\Theta(t) A_i{\bf u}_k + \Theta(1-t) A_j{\bf u}_k  \quad \hbox{for} \ \  
{\bf x}_t = t {\bf p}_i + (1-t){\bf p}_j .$$
As a consequence, if we choose any family 
$\big\{ {\bf u}_{{\bf E}} : {\bf E}$ being a mesh edge$\big\}$
of vectors
such that 
${\bf u}_{\{\bf E\}} \not\parallel {\bf E}$,
by applying the construction $(\theFS.\theFT)$ 
with the associated data over 
every triangle of the mesh,
we obtain a ${\cal C}^1$-smooth (continuously differentable)
function on
the union 
of the mesh triangles.

The shape functions $\Phi,\Theta$ and also the
RSD method described above appeared in [\MDPI]
without heuristic introduction.
Actually they arose from our earlier study 
[\STA]
with somewhat restrictive postulates
on the possible polynomial ${\cal C}^1$-interpolations
over triangular meshes
based on computer algebraic analysis of technics
developed in
[\Zlam],[\Serg],[\Hahm],[\Cao].

Our goal here
is to characterize all RSD methods
with shape functions
$\Psi_0,\Psi_1\in{\cal C}^1\big([0,1]\big)$ instead of $(\theSHS.\theSHT)$.
This requires a completely different approach as that in [\STA]
relying heavily upon polynomial identities.
From our main result 
Theorem $3.5$
%
it turns out that
there is a very general pattern 
behind $(\theFS.\theFT)$
in the form
\ADT
\newcounter{FCHIS} \setcounter{FCHIS}{\value{SEC}}
\newcounter{FCHIT} \setcounter{FCHIT}{\value{THM}}
$$\begin{aligned}
&F_0 = \sum_{k=1}^3 \Big[
\Psi_0\big(\lambda_k\big) f_k + 
\Psi_1\big(\lambda_k\big) A_k({\bf x}-{\bf p}_k) \Big] ,\cr 
&H = \hskip-3mm \sum_{
\{ \ell,m,n\}=\{1,2,3\}} 
\biggl\{ \!
f_\ell {G_\ell {\bf u}_n \over G_n {\bf u}_n}  \chi_0(\lambda_\ell,\lambda_m,\lambda_n) \ +\cr
&\hskip15mm + A_\ell({\bf p}_m\!\!-\!{\bf p}_\ell) 
{G_\ell {\bf u}_n \over G_n {\bf u}_n}
\chi_1 (\lambda_\ell,\lambda_m,\lambda_n) \!\biggr\} . 
\end{aligned}
\leqno(\NUM)
$$
If 
$\Psi_0,\Psi_1$ fulfill the minimal necessary condition
$(2.1)$ 
for giving rise to a 
Hermite interpolation of the type $F_0$ in $(\theFS.\theFT)$ of first order
then we can find plenty of suitable functions
$\chi_0,\chi_1\in{\cal C}^1_0(\RR_+^3)$,
namely those satisfying only the not too restrictive conditions 
$(3.6),\!(3.7)$,
in order for $H$ being an RSD correction for $F_0$.
This degree of freedom enables various canonical constructions for the modifiers
$\chi_0,\chi_1$ in terms of suitable shape functions $\Psi_0,\Psi_1$
to ensure specific geometric properties of the interpolation operator 
$\hbox{{\Goth F}}: \big\{ ({\bf p}_k,f_k,A_k): k=1,2,3\big\} \mapsto 
\big[ F_0+H$ by $(\theFCHIS.\theFCHIT)\big]$.  

We 
complete the paper with the investigation of two essential 
properties of the interpolation $\hbox{{\Goth F}}$
having obvious importance in applications used to construct
smooth surfaces in the form of the graph of a function $\RR^2\to\RR$
passing through a finite family of points in $\RR^3$ 
(obtained as vertices of a tringularization from scanned data):

(i) The {\it range shift property}  
$\hbox{{\Goth F}}\big\{ ({\bf p}_k,1,0) : k=1,2,3\big\} \equiv 1$, 
that is all constant functions remain invariant when 
interpolated with $\hbox{{\Goth F}}$. In this case the graph of
$\hbox{{\Goth F}}\big\{ ({\bf p}_k,f_k+c,f^\prime({\bf p}_k)) : k=1,2,3\big\}$
is just a shifted form of the graph of
$\hbox{{\Goth F}}\big\{ ({\bf p}_k,f_k,({\bf p}_k)) : k=1,2,3\big\}$
for any function $f\in{\cal C}^1(\RR^2)$.

(ii) {\it Affinity invariance}: 
$\hbox{{\Goth F}}\big\{ ({\bf p}_k,A{\bf x}+b,A) : k=1,2,3\big\} \equiv 
A{\bf x}+b$.
That is all affine (constant+linear)
functions remain invariant when 
interpolated with $\hbox{{\Goth F}}$. In particular
if the graph of a function 
${\rm Conv}\{{\bf p}_1,{\bf p}_2,{\bf p}_3\}\to\RR$ 
is a 3D plane triangle then the function obtained with interpolation by 
$\hbox{{\Goth F}}$
has the same graph. 

In Section 4 we give a parametric classification of the 
procedures with range shift property. 
Concerning affinity invariance, we have no complete results yet:
in Section 5 we provide a several necessary and sufficient algebraic conditions.   
As a by no means obvious fact, it turns out that the procedure $(\theFS.\theFT)$ with the polynomials
$(\theSHS.\theSHT)$ has range shift property but fails to be affinity invariant.
We also provide an affinity invariant procedure with polynomials of
6-th degree: namely that of the form 
$(\theFCHIS.\theFCHIT)$ with $\Psi_0 =\Psi_1 =\Phi$ and 
$\chi_0(t_1,t_2,t_3)=30 t_1^2 t_2^2 t_3$,
$\chi_1(t_1,t_2,t_3)=30 t_1^2 t_2^3 t_3 + 15 t_1^2 t_2^2 t_3^2$.

\vskip6mm
\ADS
\noindent{\bfs \theSEC. Preliminaries}
\medskip
\rm

For standard terminology,
$\RR^N$ resp. $\RR^N_+$ stand for the set of all 
real resp. non-negative $N$-tuples 
$[\xi_1,\ldots,\xi_N]$. 
Actually, we shall only be interested in the cases of 
dimensions $N=1,2,3$.

Given a subset $\Omega$ of $\RR^N$, 
we write ${\cal C}(\Omega)$ for the family of all
continuous functions with domain $\Omega$.
If the interior $\Omega^{\circ}$ of $\Omega$ is dense in $\Omega$,
we define ${\cal C}^1(\Omega)$ as the set 
consisting of all functions $f \in {\cal C}(\Omega)$
with continuous partial derivatives $D_1 f, \ldots,D_N f$
on $\Omega^\circ$ {\it admitting a continuous extension to
$\Omega$} with the value denoted also by 
$D_k f({\bf x})$ at the points ${\bf x}\in\Omega\setminus\Omega^\circ$.
As a folklore consequence of Whitney's extension theorem [\WHIT],
if $\Omega$ is closed in $\RR^N$ then
every function $f\in{\cal C}^1(\Omega)$ can be regarded as the
restriction of a continuously differentiable function
defined on $\RR^N$.

For any function $f\in{\cal C}^1(\Omega)$, we write \
$f^\prime$ \  
for its {\it Fr\'echet derivative}
defined at any point ${\bf x}\in\Omega$ 
as
the linear functional
\vskip-5mm 
$$f^\prime({\bf x}) {\bf u} =
\sum_{k=1}^N D_k f({\bf x}) \upsilon_k
\qquad 
\big( {\bf u}=[\upsilon_1,\ldots,\upsilon_N] \big).$$
In particular 
$f^\prime({\bf x}){\bf u} = 
\lim\limits_{t\to 0} t^{-1} \big[ f({\bf x}+t{\bf u}) -f{\bf u}\big]$ 
is the familiar 
directional derivative     
at the points ${\bf x}\in\Omega^\circ$.
If $\Omega$ is a closed set, we define
${\cal C}_0^k(\Omega) = \big\{ f\in{\cal C}^k(\Omega) :
f({\bf x})=0 \ \big({\bf x}\in\Omega\setminus\Omega^o\big)\big\}$ $(k=0,1)$.

Given any subset ${\bf P} \subset \RR^2$,
we write 
${\rm Conv}({\bf P})$ 
for its convex hull. 
If ${\bf P} =\{ {\bf p}_1,{\bf p}_2,{\bf p}_3\}$
and the triangle   
${\bf T}={\rm Conv}({\bf P})$ 
is non-degenerate, 
the {\it barycentric weights} (weight functions
associated to the vertices)
[\Berger]
are the functions 
\vskip-5mm
$$\lambda_i=\lambda_{{\bf p}_i}^{{\bf T}} : {\bf x} \mapsto  
{ {\rm det}\big[ {\bf x} - {\bf p}_j, {\bf x} - {\bf p}_k\big] 
\over
{\rm det}\big[ {\bf p}_i - {\bf p}_j, {\bf p}_i - {\bf p}_k\big] }
\quad \hbox{with} \ (i,j,k) \in S_3$$
\vskip-7mm
\noindent 
where
\vskip-3mm  
$$S_3 =
\big\{ (1,2,3),(2,1,3),(2,3,1),(3,2,1),(3,1,2),(1,3,2)\big\}$$
denotes the set of all pemutations of the indices $1,2,3$.
For the sake of brevity, in the sequel we omit the background parameters
in most formulas
(like writing $\lambda_i \!=\! \lambda_{{\bf p}_i}^{{\bf T}}$ above) 
without danger of confusion.
It is well-known that
the mapping 
${\bf x} \mapsto [\lambda_1({\bf x}),\lambda_2({\bf x}),\lambda_3({\bf x})]$
is a homeomorphism between ${\bf T}$ and the {\it 3D-unit simplex}
\vskip-5mm
$$\Delta_3 := \big\{ [t_1,t_2,t_3] \in\RR_+^3 : t_1+t_2+t_3 =1 \big\} ,$$
moreover
$[\lambda_1({\bf x}),\lambda_2({\bf x}),\lambda_3({\bf x})]$
is the unique triple of non-negative numbers such that 
${\bf x}\!= \!\sum_k \lambda_k({\bf x}) {\bf p}_k$ and 
$\sum_k  \lambda_k({\bf x}) \!=\! 1$.
Furthermore 
$\lambda_1,\lambda_2,\lambda_3:\RR^2\!\to\!\RR$
are affine (linear+constant) functions with the necessarily constant 
Fr\'echet derivatives 
$$G_i {\bf u} :=\lambda_i^\prime ({\bf x}) = 
{d\over dt}\Big\vert_{t=0} \lambda_i({\bf x}+ t{\bf u}) \quad
\hbox{independently of } {\bf x}.$$
Given any pair 
${\bm \Psi} = [\Psi_0,\Psi_1]$ with
$\Psi_0,\Psi_1 \in{\cal C}^1\big([0,1])$
such that
\ADT
\newcounter{PSIS} \setcounter{PSIS}{\value{SEC}}
\newcounter{PSIT} \setcounter{PSIT}{\value{THM}}
$$\Psi_0(0) = \Psi_0^\prime(0) = \Psi_1(0) = 
\Psi_1^\prime(0)= \Psi_0^\prime(1) = 0, \quad
\Psi_0(1) =\Psi_1(1)=1 
\leqno(\NUM)$$  
we introduce the 
{\it basic triangular interpolation of first order
with the shape functions $\Psi_0,\Psi_1$}
as the operator
\ADT
\newcounter{FOS} \setcounter{FOS}{\value{SEC}}
\newcounter{FOT} \setcounter{FOT}{\value{THM}}
$$\FFg_0^{\bm \Psi} :
\big\{ ({\bf p}_k,f_k,A_k) \big\}_{k=1}^3 \mapsto
\sum_{k=1}^3 \Big\{ \Psi_0\big(\lambda_k\big) f_k \!+\! \Psi_1\big(\lambda_k\big) A_k\big({\bf x} \!-\!{\bf p}_k\big) \Big\} 
\leqno(\NUM)$$
defined for all {\it first order function germs}
$\big\{ ({\bf p}_k,f_k,A_k) : k \!=\!1,2,3\big\}$
with ${\bf T} ={\rm Conv}\{ {\bf p}_1, {\bf p}_2, {\bf p}_1\} \subset\RR^2$
being a non-degenerate triangle, 
$f_k\in\RR$ and $A_k\in {\cal L}(\RR^2,\RR)$ in terms of the weights
$\lambda_k = \lambda_k^{\bf T}$.  
By definition, the domain of the function 
$F_0 = \FFg_0^{\bm \Psi} \big\{ ({\bf p}_k,f_k,A_k) \big\}_{k=1}^3$ 
is only the triangle ${\bf T}$ and $F_0 \in {\cal C}^1({\bf T})$.
By straightforward calculation,
for its Fr\'echet derivative 
we have
\ADT
\newcounter{DERS} \setcounter{DERS}{\value{SEC}}
\newcounter{DERT} \setcounter{DERT}{\value{THM}}
$$F_0^\prime =\sum_{k=1}^3 \Big\{ 
\big[\Psi_0^\prime\big(\lambda_k\big)  f_k + 
\big[\Psi_1^\prime\big(\lambda_k\big) A_k\big({\bf x}-{\bf p}_k\big) \big] G_k +
\Psi_1\big(\lambda_k\big) A_k \Big\} .
\leqno(\NUM)$$

\ADT
\newcounter{REMiS} \setcounter{REMiS}{\value{SEC}}
\newcounter{REMiT} \setcounter{REMiT}{\value{THM}}
\bigskip
\noindent
{\bf \NUM. Remark.} (a) Conditions (\thePSIS.\thePSIT) 
do not restrict the value of $\Psi_1^\prime(1)$.

(b) In view of $(\theDERS.\theDERT)$, it is not hard to see that 
$(\thePSIS.\thePSIT)$ is suffiecient and necessary for 
$F_0$   
to satisfy the relations
\ADT
\newcounter{POSTS} \setcounter{POSTS}{\value{SEC}}
\newcounter{POSTT} \setcounter{POSTT}{\value{THM}}
$$F_0\big({\bf p}_i \big) =f_i, \quad F_0^\prime\big({\bf p}_i\big) = A_i 
\qquad (i=1,2,3) 
\leqno(\NUM)$$
under any choice of 
the first order function germs $\big\{ {\bf p}_k,f_k,A_k) \big\}_{k=1}^3$.

(c) We can express the terms $A_k({\bf x}-{\bf p}_k)$ appearing in 
$(\theFCHIS.\theFCHIT)$, $(\theFOS.\theFOT)$ resp. $(\thePOSTS.\thePOSTT)$ 
in the form of linear combination of the weights due to the identity
${\bf x} = \sum_{i=1}^3 \lambda_i({\bf x}) {\bf p}_i$ as
\vskip-4mm 
$$A_k({\bf x}-{\bf p}_k) = \sum_{i=1}^3 A_k({\bf p}_i - {\bf p}_k) \lambda_i .$$ 

\ADT
\newcounter{DEFiS} \setcounter{DEFiS}{\value{SEC}}
\newcounter{DEFiT} \setcounter{DEFiT}{\value{THM}}
\noindent
{\bf \NUM. Definition.} Henceforth we say that 
${\bm \Psi} =[\Psi_0,\Psi_1]$ 
is a {\it pair of admissible shape functions} if 
$\Psi_0,\Psi_1\in{\cal C}^1([0,1])$ 
and 
$(\thePSIS.\thePSIT)$ holds. 
We also say for short that $\big\{ ({\bf p}_k,f_k,A_k) \big\}_{k=1}^3$
is an {\it admissible  function germ} if
${\bf p}_1,{\bf p}_2,{\bf p}_3 \in\RR^2$ form a non-degenete triangle,
$f_k\in\RR$ and $A_k\in{\cal L}(\RR^2,\RR)$.

\bigskip

\noindent
As an immediate consequence of $(\theDERS.\theDERT)$, since
\ADT
\newcounter{PARS} \setcounter{PARS}{\value{SEC}}
\newcounter{PART} \setcounter{PART}{\value{THM}} 
$$\begin{aligned}
&\lambda_i\big( {\bf p}_j\big) =\delta_{ij}, \quad
G_i \big({\bf p}_j -{\bf p}_k \big) = \delta_{ij} - \delta_{ik} ,\cr 
&\lambda_i({\bf x}_t^k)=t, \quad \lambda_j({\bf x}_t^k)=1-t, \quad 
\lambda_k({\bf x}_t^k)=0, \cr 
&{\bf x}^k_t-{\bf p}_i = (1-t) ({\bf p}_j-{\bf p}_i), \quad 
{\bf x}^k_t-{\bf p}_j = t ({\bf p}_i-{\bf p}_j) , 
\end{aligned}
\leqno(\NUM)$$ 
we obtain the following observation. 

\bigskip
\ADT
\newcounter{LEMiS} \setcounter{LEMiS}{\value{SEC}}
\newcounter{LEMiT} \setcounter{LEMiT}{\value{THM}}
\noindent
{\bf \NUM. Lemma.} \it Given an admissible pair 
${\bm \Psi} = [\Psi_0,\Psi_1]$
of shape functions along with an admissible function germ 
$\hbox{\Goth g}=\big\{ ({\bf p}_k,f_k,A_k) \big\}_{k=1}^3$,
if $i,j,k\in\{ 1,2,3\}$ are three different indices
then, 
at the generic point
\ADT
\newcounter{GENPS} \setcounter{GENPS}{\value{SEC}}
\newcounter{GENPT} \setcounter{GENPT}{\value{THM}}
$${\bf x}_t^k := t{\bf p}_i + (1-t){\bf p}_j \qquad 
(0\le t\le 1)
\leqno(\NUM)$$
of the edge $[{\bf p}_i,{\bf p}_j]$ in the triangle 
${\bf T} ={\rm Conv}\{ {\bf p}_1,{\bf p}_2,{\bf p}_3\}$, 
for the function
$F_0 =\FFg^{\bm \Psi}_0 \hbox{\Goth g}$  
we have
\ADT
\newcounter{FnulS} \setcounter{FnulS}{\value{SEC}}
\newcounter{FnulT} \setcounter{FnulT}{\value{THM}} 
$$\begin{aligned}
F_0\big({\bf x}_t^k\big) &= 
\Psi_0(t) f_i + \Psi_1(t) (1-t) A_i ({\bf p}_j - {\bf p}_i) +\cr 
&\hskip5mm +  
\Psi_0(1-t) f_j + \Psi_1(1-t) t A_j ({\bf p}_i - {\bf p}_j) ; \cr
F_0^\prime\big({\bf x}_t^k\big) &=  
\big[ \Psi_0^\prime(t) f_i + 
\Psi_1^\prime(t) (1-t)A_i({\bf p}_j-{\bf p}_i)\big] G_i
+ \Psi_1(t) A_i + \cr 
&\hskip-7mm + \big[ \Psi_0^\prime(1-t) f_j + 
\Psi_1^\prime(1-t) t A_j({\bf p}_i-{\bf p}_j)\big] G_j
+ \Psi_1(1-t) A_j .
\end{aligned}
\leqno(\NUM)$$
\rm

\vskip10mm
\ADS
\noindent{\bfs \theSEC. 
 Generic algorithm of reduced side derivatives (RSD)}
\bigskip
\rm

In this section we are looking for
${\cal C}^1$-spline interpolation procedures 
analogous to those described in [\MDPI] but
with more general shape functions 
$\Psi_0,\Psi_1 \in{\cal C}^1[0,1]$ instead of   
$\Phi,\Theta$ there.

Henceforth 
let ${\bf p}_1,{\bf p}_2,{\bf p}_3 \in \RR^2$ be the vertices
of some (arbitrarily fixed)
non-degenerate triangle ${\bf T}$ with respective weight functions
$\lambda_m:\RR^2\!\to\!\RR$ and derivative weights 
$G_m \!:= \lambda_m^\prime \!\in\! {\cal L}(\RR^2,\RR)$ $(m\!=\! 1,2,3)$  
and let ${\bf u}_1,{\bf u}_2,{\bf u}_3 \in \RR^2$ 
be given vectors such that
${\bf u}_k \not\parallel {\bf p}_i-{\bf p}_j$ whenever $(i,j,k)\in S_3$.

Furthermore
${\bm \Psi} \!=\! [\Psi_0,\Psi_1]$ resp.
$\hbox{\Goth g} \!=\!\big\{ ({\bf p}_k,\! f_k,\! A_k) \big\}_{k=1}^3$ 
denote a 
fixed admissible pair of shape functions
and  
a function germ.
As earlier, we write\break  
\vskip-7mm
$$F_0 = \hbox{\Goth F}_0^{\bm \Psi} \hbox{\Goth g} .$$
Our starting point is the following immediate consequence
of $(\theFnulS.\theFnulT)$:  

\ADT
\newcounter{REMiiS} \setcounter{REMiiS}{\value{SEC}}
\newcounter{REMiiT} \setcounter{REMiiT}{\value{THM}}
\bigskip
\noindent
{\bf \NUM. Remark.} 
Let $H\in{\cal C}^1_0({\bf T})$.
Then 
the function $F:=F_0-H$ coincides with $F_0$ along the edges of ${\bf T}$
and hence, at the vertices, it 
has also the properties $F({\bf p}_i) =f_i$, $F({\bf p}_i) =A_i$
$(i=1,2,3)$ analogous to $(\thePOSTS.\thePOSTT)$.

\ADT
\newcounter{DEFiiS} \setcounter{DEFiiS}{\value{SEC}}
\newcounter{DEFiiT} \setcounter{DEFiiT}{\value{THM}}
\bigskip
\noindent
{\bf \NUM. Definition.}  
We say that the function
\ADT
\newcounter{MODS} \setcounter{MODS}{\value{SEC}}
\newcounter{MODT} \setcounter{MODT}{\value{THM}} 
$$F = F_0- H , \quad  H\in{\cal C}^1_0({\bf T}) 
\leqno(\NUM)$$ 
is an 
{\it RSD modification} of $F_0$ (with respect to the directions 
${\bf u}_1,{\bf u}_2,{\bf u}_3$ along the edges 
$\big[{\bf p}_2,{\bf p}_3\big], \big[{\bf p}_3,{\bf p}_1\big],
\big[{\bf p}_1,{\bf p}_2\big]$ respectively) 
if
for any index pemutation $(i,j,k)\in S_3$, 
$$F^\prime \big({\bf x} \big) {\bf u}_k \!\!=\!
\Psi_1\big(\lambda_i({\bf x})\big) A_i {\bf u}_k + 
\Psi_1\big(\lambda_j({\bf x})\big) A_j {\bf u}_k \ 
\hbox{whenever} \ 
{\bf x}\!\in\!\! \big[{\bf p}_i,{\bf p}_j\big] .$$
Conveniently, in this case we refer to $H$ as an 
{\it RSD modifier of $F_0$}.
According to Lemma $\theLEMiS.\theLEMiT$, we have the following.

\ADT
\newcounter{CORiS} \setcounter{CORiS}{\value{SEC}}
\newcounter{CORiT} \setcounter{CORiT}{\value{THM}}
\bigskip
\noindent
{\bf \NUM. Corollary.} {\it
A function $H \in{\cal C}^1_0({\bf T})$ is an 
RSD modifier of $F_0$ 
with respect to the directions 
${\bf u}_k$ along the edges $\big[{\bf p}_i,{\bf p}_j\big]$
if and only if
$$\begin{aligned}
H^\prime\big({\bf x}_t^k\big) {\bf u}_k &=  
\big[ \Psi_0^\prime(t) f_i + 
\Psi_1^\prime(t) (1-t)A_i({\bf p}_j-{\bf p}_i)\big] [G_i {\bf u}_k] + \cr 
&+ \big[ \Psi_0^\prime(1-t) f_j + 
\Psi_1^\prime(1-t) t A_j({\bf p}_i-{\bf p}_j)\big] [G_j {\bf u}_k] 
\end{aligned}$$
whenever ${\bf x}_t^k = t {\bf p}_i + (1-t){\bf p_j}$ with
$(i,j,k)\in S_3$ and $0\le t\le 1$.
}

\bigskip
\ADT
\newcounter{THMiS} \setcounter{THMiS}{\value{SEC}}
\newcounter{THMiT} \setcounter{THMiT}{\value{THM}}
\noindent
{\bf \NUM. 
{\bf Theorem.} {\it
Let $\Omega \subset\RR^3$ be a set 
whose interior contains
$\Delta_3\cap (0,1)^3$.
Assume 
$\chi_0,\chi_1 \in{\cal C}^1(\Omega)$  
are functions vanishing along the edges
$\Delta_{3,k} \!=\! \big\{ (t_1,t_2,t_3) \!\in\!\Delta_3: t_k \!=\! 0\big\}$
$(k=1,2,3)$ of $\Delta_3$
such that
\ADT
\newcounter{CHIS} \setcounter{CHIS}{\value{SEC}}
\newcounter{CHIT} \setcounter{CHIT}{\value{THM}}
$$D_3\chi_0(t,1-t,0)=\Psi_0^\prime(t), \quad 
D_3\chi_1(t,1-t,0) = \Psi_1^\prime(t)\cdot(1-t)
\leqno(\NUM)$$
with the following marginal conditions on the derivatives 
\ADT
\newcounter{CHImS} \setcounter{CHImS}{\value{SEC}}
\newcounter{CHImT} \setcounter{CHImT}{\value{THM}}
$$\chi_r^\prime({\bf t})=0 \ \ 
\big({\bf t}\in\Delta_{3,1}\cup\Delta_{3,2}\big) ; \quad
D_m \chi_r({\bf t})=0 \ \ \big({\bf t}\in\Delta_{3,3}, \ 
m=1,2\big) .
\leqno(\NUM)$$
Then the function $H$ in $(\theFCHIS.\theFCHIT)$
is an RSD modifier for $F_0$ 
under any choice of the vectors ${\bf u}_1,{\bf u}_2,{\bf u}_3$
with \ ${\bf u}_i \not\parallel {\bf p}_j-{\bf p}_k$
$\big( (i,j,k)\in S_3 \big)$.
}

\bigskip
{\bf Proof.} \rm
Consider any permutation $(i,j,k)\in S_3$.
According to Corollary $\theCORiS.\theCORiT$, it suffices to see that,
at the points
${\bf x}_t = t{\bf p}_i + (1-t){\bf p}_j$ 
we have
\ADT
\newcounter{RSDmS} \setcounter{RSDmS}{\value{SEC}}
\newcounter{RSDmT} \setcounter{RSDmT}{\value{THM}}
$$\begin{aligned}
H({\bf x}_t) &= 0 ,\cr  
H^\prime({\bf x}_t) {\bf u}_k &= 
\big[ \Psi_0^\prime(t) f_i + 
\Psi_1^\prime(t) (1-t)A_i({\bf p}_j-{\bf p}_i)\big] G_i {\bf u}_k + \cr 
&\hskip6mm + \big[ \Psi_0^\prime(1-t) f_j + 
\Psi_1^\prime(1-t) t A_j({\bf p}_i-{\bf p}_j)\big] G_j {\bf u}_k .
\end{aligned}
\leqno(\NUM)$$
Notice that, in terms of the triples
$${\bf t}_{\ell,m,n} \!=\!
\big( \lambda_\ell({\bf x}_t),\lambda_m({\bf x}_t),\lambda_n({\bf x}_t)\big)
\quad \big( (\ell,m,n) \in S_3\big) ,$$
and by setting 
$(\ell,m,n)_1 \!=\!\ell,\, (\ell,m,n)_2 \!=\! m,\, (\ell,m,n)_3 \!=\! n$,
we can write
$$\begin{aligned}
&H({\bf x}_t) = \sum_{(\ell,m,n)\in S_3} 
{G_\ell {\bf u}_n \over G_n {\bf u}_n} 
\bigg[ f_\ell \chi_0({\bf t}_{\ell,m,n}) +
A_\ell ({\bf p}_m - {\bf p}_\ell) \chi_1({\bf t}_{\ell,m,n}) \bigg] ,\cr 
&H^\prime({\bf x}_t) {\bf u}_k  = \hskip-2mm 
\sum_{(\ell,m,n) \in S_3} \hskip-2mm {G_\ell {\bf u}_n \over G_n {\bf u}_n} 
\bigg[ f_\ell \Big[ \sum_{q=1}^3 D_q \chi_0({\bf t}_{\ell,m,n}) 
G_{(\ell,m,n)_q} {\bf u_k}  \Big] +\cr 
&\hskip32mm + 
A_\ell({\bf p}_m-{\bf p}_\ell) \Big[
\sum_{q=1}^3 D_q \chi_1({\bf t}_{\ell m n}) 
G_{(\ell,m,n)_q} {\bf u_k} \Big] \bigg\} .
\end{aligned}$$
Observe that we have
${\bf t}_{ijk} \!=\! (t,1-t,0) \!\in\! \Delta_{3,3}$, 
${\bf t}_{jik} \!=\! (1-t,t,0) \!\in\! \Delta_{3,3}$,
${\bf t}_{ikj} \!=\! (t,0,1-t) \!\in\! \Delta_{3,2}$,
${\bf t}_{kij} \!=\! (0,t,1-t) \!\in\! \Delta_{3,2}$,
${\bf t}_{jki} \!=\! (1-t,0,t) \!\in\! \Delta_{3,1}$,
${\bf t}_{kji} \!=\! (0,1-t,t) \!\in\! \Delta_{3,1}$.
Hence the relation $H({\bf x}_t)=0$ is immediate.  
According to $(\theCHImS,\theCHImT)$,
all the terms $D_q \chi_r({\bf t}_{\ell,m,n})$ vanish 
except for the cases with
$\lambda_n({\bf x}_t) =0$ and $q=3$ that is for
$(q,r,\ell,m,n) \in\big\{ (3,r,i,j,k),(3,r,j,i,k): r=0,1\big\}$.
In view of $(\theCHIS.\theCHIT)$ we have
$D_3\chi_0({\bf t}_{ijk}) =\Psi_0^\prime(t)$,
$D_3\chi_0({\bf t}_{jik}) =\Psi_0^\prime(1-t)$,
$D_3\chi_1({\bf t}_{ijk}) =\Psi_1^\prime(t)(1-t)$,
$D_3\chi_1({\bf t}_{jik}) =\Psi_1^\prime(1-t)t$
which completes the proof of $(\theRSDmS.\theRSDmT)$
and hence the theorem.

\bigskip

\ADT
\newcounter{DEFiiiS} \setcounter{DEFiiiS}{\value{SEC}}
\newcounter{DEFiiiT} \setcounter{DEFiiiT}{\value{THM}}
\noindent
{\bf \NUM. Definition.} We shall say that
$[\Psi_0,\Psi_1,\chi_0,\chi_1]$ is an {\it RSD tuple} if
$[\Psi_0,\Psi_1]$ is an admissible pair of shape functions
and 
$\chi_0,\chi_1$ are functions
satisfying the requirements of Theorem $\theTHMiS.\theTHMiT$. 

\DEL{
\bigskip

\ADT
\newcounter{REMxivS} \setcounter{REMxivS}{\value{SEC}}
\newcounter{REMxivT} \setcounter{REMxivT}{\value{THM}}
\noindent
{\bf \NUM. Remark.} The (numerically less interesting)
converse of Theorem $\theTHMiS.\theTHMiT$
is also valid:

{\it Given a 
non-degenerate triangle ${\rm Conv}\{ {\bf p}_1, {\bf p}_2, {\bf p}_3\}$,
a pair $\Psi_0,\Psi_1$
of admissible shape functions along with a pair 
$\chi_0,\chi_1\in{\cal C}^1(\Omega)$ of functions
defined on a domain 
$\Omega\subset\RR^3$ such that 
${\rm Closure}(\Omega^o) \supset \Delta_3$,
if the function $H$ is an RSD modifier of $F_0$
under any choice of $f_1,f_2,f_3\in\RR$, $A_1,A_2,A_3\in{\cal L}(\RR^2,\RR)$
and ${\bf u}_1,{\bf u}_2,{\bf u}_3\in\RR^2$ with 
${\bf u}_k\not\parallel {\bf p}_j-{\bf p}_i$ $\big( (i,j,k)\in S_3\big)$
in $(\theFCHIS.\theFCHIT)$, 
then we necessarily have
$(\theCHIS.\theCHIT),(\theCHImS.\theCHImT)$.
}

The proof relies upon a careful study of the second
equation in $(\theRSDmS.\theRSDmT)$
under suitable choices of 
$f_\ell\in\{ 0,1\}$, $A_\ell\in\{ 0, G_i : i=1,2,3\}$
such that only a unique term does not vanish from 
$\big[ A_\ell({\bf p}_m-{\bf p}_\ell) : \ell,m=1,2,3\big]$ 
and taking ${\bf u}_n$ $(n=1,2,3)$ in the form
${\bf u}_n^{\omega} = [{\bf p}_n-({\bf p}_\ell+{\bf p}_m)/2]+
\omega [{\bf p}_m-{\bf p}_\ell]$ for 
$(\ell,m,n)=(1,2,3),2,3,1),(3,1,2)$.

We assume without loss of generaliy that $i=1,\, j=2,\, k=3$
and we let
${\bf u}_1={\bf u}_1^\alpha$, ${\bf u}_1={\bf u}_1^\beta$, 
${\bf u}_1={\bf u}_1^\gamma$
with $G_n{\bf u}_n=1$, $G_1{\bf u}_2=-1/2+\beta$, $G_1{\bf u}_3=-1/2-\gamma$,
$G_2{\bf u}_3=-1/2+\gamma$, $G_2{\bf u}_1=-1/2-\alpha$,
$G_3{\bf u}_1=-1/2+\alpha$, $G_3{\bf u}_2=-1/2-\beta$.

With the choice $f_1=1$, $f_2=f_3=0$, $A_\ell=0$, 
we deduce that
$$\begin{aligned}
&\big(\!\! -\!{1\over 2}\!-\! \gamma\!\big) \Big\{  
D_1\chi_0({\bf t}_{1,2,3})\big(\!\! -\!{1\over 2}\!-\!\gamma \!\big) \!+\!
D_2\chi_0({\bf t}_{1,2,3})\big(\! -{1\over 2}\!+\! \gamma\big) \!+\!
D_3\chi_0({\bf t}_{1,2,3}) \Big\} +\cr
&+ \big(\!\! -\!{1\over 2}\!+\!\beta\big) \Big\{  
D_1\chi_0({\bf t}_{1,3,2})\big(\!\! -{1\over 2}\!-\!\gamma\big) \!+\!
D_2\chi_0({\bf t}_{1,3,2}) \!+\! D_3\chi_0({\bf t}_{1,3,2}) \Big\} \!=\!\cr
&=\Psi_0^\prime(t)\big(\! -{1\over 2}\!+\!\gamma \!\big) 
\end{aligned}
$$
}

\bigskip

\ADT
\newcounter{LEMiiS} \setcounter{LEMiiS}{\value{SEC}}
\newcounter{LEMiiT} \setcounter{LEMiiT}{\value{THM}}
\DEL{
\noindent
{\bf \NUM. Lemma.} {\it 
Given any function $h \!\in\! {\cal C}_0(\RR^2_+)$,
there exists
$\chi \!\in\! {\cal C}^1_0(\Omega)$ on 
$\Omega=\RR_+^3\setminus\{ (0,0,0)\}$ 
such that for all $t_1,t_2,t_3\ge 0$ we have  
$$\begin{aligned}
&0 = D_m \chi(0,t_2,t_3)=D_m(t_1,0,t_3) = D_m(t_1,t_2,0) \quad (m=1,2), \cr
&h(t_1,t_2) = D_3\chi(t_2,t_2,0), \ \ 
0 = D_3\chi(0,t_2,t_3) = D_3\chi(t_1,0,t_3) .
\end{aligned}$$
} 
\hskip5mm {\bf Proof.} 
For ${\bf t} = (t_1,t_2,t_3) \in\RR^3_+$,
define the smoothing 
$$\widehat{h}({\bf t}) \!=\! {1\over t_3^2} \!\int_{s_1=t_1}^{t_1+t_3} \!
\int_{s_2=t_2}^{t_2+t_3} \!\! h(s_1,s_2) \, ds_2 ds_1  \ \ (t_3 \!>\! 0), \quad  \widehat{h}({\bf t}) \!=\! h(t_1,t_2)
 \ \ (t_3 \!=\! 0)$$
of $h(t_1,t_2,0)$ and fix any function \ $\Phi\in{\cal C}^1\big( [0,1]\big)$ \ 
such that
$$\phi(t) =0 \ \ \big(0\le t\le {1\over 3}\big), \quad 
\phi^\prime(t) \ge 0 \ \ \big({1\over 3}\le t\le {2\over 3}\big), \quad
\phi(t) =0 \ \ \big({2\over 3}\le t\le 1\big) .$$
We show that the function
$$\chi(t_1,t_2,t_3) = 
\phi\Big({t_1\over t_1+t_3}\Big) \phi\Big({t_2\over t_2+t_3}\Big)
t_3\widehat{h} (t_1,t_2,t_3)$$
extended with $\chi(t_1,t_2,0)=0$ suits the requirements of the lemma.

It is folklore that the function $\widehat{h}$ is continuous, and
by the Newton-Leibniz theorem,
it is continuously differentiable
at the points 
of the wedge region 
$\RR^2_+\times\RR_{++} \big(=
\{ (t_1,t_2,t_3): t_1,t_2\ge 0, \ t_3>0 \} \big)$.
By construction, on the relatively open subdomains
$$\begin{aligned}
&\Omega_0 = \Big\{ (t_1,t_2,t_3)\in\Omega : 
\min \Big({t_1\over t_1+t_3}, {t_2\over t_2+t_3}\Big) < {1\over 3}\Big\} ,\cr 
&\Omega_1 = \Big\{ (t_1,t_2,t_3)\in\Omega : 
\min \Big({t_1\over t_1+t_3}, {t_2\over t_2+t_3}\Big) > {2\over 3}\Big\}
\end{aligned}$$  
we have 
$$\chi({\bf t}) =0 \quad ({\bf t}\in\Omega_0), \quad
\chi({\bf t}) = t_3\widehat{h}({\bf t}) \quad ({\bf t}\in\Omega_1) .$$
In particular \
$\chi\in{\cal C}_0(\Omega)$ and the restriction of $\chi$ to 
$\RR^2_+\times\RR_{++}$ is continuosly differentiable with 
$D_m\chi({\bf t})=0$ for ${\bf t}\in \RR^2_+\times\RR_{++}$
if $t_1=0$ or $t_2=0$.
Therefore it only remains to prove that
given any point ${\bf t}^*=(0,t_2^*,t_3^*)$ with 
$t_1^*,t_2^*\ge 0$ and $t_1^*+t_2^*>0$,
we have
$$D_m [t_3\widehat{h}]({\bf t}) \!\!\to\! 0 \ (m \!=\!\! 1,2) ,\ 
D_3 [t_3\widehat{h}]({\bf t}) \!\!\to\! h(t_1^*,t_2^*)
\ \ \hbox{whenever} \  
\RR^2_+\!\times\RR_{++} \!\ni\! {\bf t} \!\to\! {\bf t}^* .$$
However, this fact is an immedite consequence of the mean value expressions
$$\begin{aligned}
{\partial [t_3\widehat{h}]\over \partial t_1} &=
{1 \over t_3} \int_{s_2=t_2}^{t_2+t_3}
\big[ h(t_1+t_3,s_2) - h(t_1,s_2) \big] \, ds_2 =\cr 
&= h\big( t_1 + t_3,r_2({\bf t})\big) - h\big( t_1,q_2({\bf t})\big), \cr
{\partial [t_3\widehat{h}]\over \partial t_2} &=
{1 \over t_3} \int_{s_1=t_1}^{t_1+t_3}
\big[ h(s_1,t_2+t_3) - h(s_1,t_2) \big] \, ds_1 =\cr 
&= h\big( r_1({\bf t}),t_2+t_3\big) - h\big( q_1({\bf t}),t_2\big), \cr
{\partial [t_3\widehat{h}]\over \partial t_3} &=
-{1 \over t_3^2} \int_{s_1=t_1}^{t_1+t_3} \int_{s_2=t_2}^{t_2+t_3}
h(s_1,s_2) \ d s_2 \ d s_1 \ +\cr
&\hskip5mm + {1 \over t_3} \int_{s_2=t_2}^{t_2+t_3}
h(t_1+t_3,s_2) \ d s_2 \ +
{1 \over t_3} \int_{\! s_1=t_1}^{t_1+t_3} h(s_1,t_2+t_3) \ d s_1 =\cr 
&= -h\big( p_1({\bf t}), p_2({\bf t})\big) +
h\big( t_1+t_3,r_2({\bf t})\big) +
h\big( r_1({\bf t}),t_2+t_3\big)
\end{aligned}$$
with suitable 
$$p_1({\bf t}),q_1({\bf t}), r_1({\bf t}) \in [t_1,t_1+t_3], \quad 
p_2({\bf t}),q_2({\bf t}), r_2({\bf t}) \in [t_2,t_2+t_3] .$$
}

\DEL{
\newcounter{LEMiiS} \setcounter{LEMiiS}{\value{SEC}}
\newcounter{LEMiiT} \setcounter{LEMiiT}{\value{THM}}
\noindent
}

{\bf \NUM. Lemma.} {\it 
Given any function $h \!\in\! {\cal C}_0(\RR^2_+)$,
there exists
$\chi \!\in\! {\cal C}^1_0(\RR_+^3)$ 
such that for all $t_1,t_2,t_3\ge 0$ we have  
$$\begin{aligned}
&0 = D_m \chi(0,t_2,t_3)=D_m(t_1,0,t_3) = D_m(t_1,t_2,0) \quad (m=1,2), \cr
&h(t_1,t_2) = D_3\chi(t_2,t_2,0), \ \ 
0 = D_3\chi(0,t_2,t_3) = D_3\chi(t_1,0,t_3) .
\end{aligned}$$
} 
\hskip5mm {\bf Proof.} 
For ${\bf t} = (t_1,t_2,t_3) \in\RR^3_+$,
define the smoothing 
$$\widehat{h}({\bf t}) \!=\! {1\over t_3^2} \!\int_{s_1=t_1}^{t_1+t_3} \!
\int_{s_2=t_2}^{t_2+t_3} \!\! h(s_1,s_2) \, ds_2 ds_1  \ \ (t_3 \!>\! 0), \quad  \widehat{h}({\bf t}) \!=\! h(t_1,t_2)
 \ \ (t_3 \!=\! 0)$$
of $h(t_1,t_2,0)$ and let us fix a function
$\phi\in {\cal C}^1\big( [0,1]\big)$ 
(e.g. $\phi(t) = 3t^2-2t^3$)
with bounded derivative and being such that
$\phi(0)=\phi^\prime(0)=\phi^\prime(1)=0$ and $\phi(1)=1$.
We show that the function
$$\chi(t_1,t_2,t_3) = 
\phi\Big({t_1\over t_1+t_3}\Big) \phi\Big({t_2\over t_2+t_3}\Big)
t_3\widehat{h} (t_1,t_2,t_3) \qquad
(t_1,t_2\ge 0,\ t_3>0)$$
extended with $\chi(t_1,t_2,0)=0$ suits the requirements of the lemma.

It is folklore that,
by the Newton-Leibniz theorem,
the function $\widehat{h}$ is continuous,
moreover 
its restriction to
$\RR^2_+\times\RR_{++} \big(=
\{ (t_1,t_2,t_3): t_1,t_2\ge 0, \ t_3>0 \} \big)$
is ${\cal C}^1$ smooth. 
In particular, for any point
${\bf t}=(t_1,t_2,t_3)$ with $t_3>0$ and $t_1=0$ or $t_2=0$
we have 
$\chi({\bf t}) =0$ and $\chi^\prime({\bf t})=0$.
Furthermore, for any point ${\bf t}\in\RR_+\times\RR_{++}$, 
with the indices $m=1,2$ resp. 3
we can write
$$\begin{aligned}
&D_m\chi({\bf t}) = \!\Phi\Big({t_{3-m}\over t_{3-m} \!+\! t_3}\Big) \bigg[
\Phi^\prime\Big({t_m\over t_m \!+\! t_3}\Big){t_3^2 \over(t_m \!+\! t_3)^2}
\widehat{h}({\bf t}) \!+\!
\Phi\Big({t_m\over t_m \!+\! t_3}\Big) 
{\partial [t_3\widehat{h}] \over \partial t_1} \bigg],\cr 
&D_3\chi({\bf t}) =
-\Phi^\prime\Big({t_1 \over t_1+t_3}\Big) 
\Phi\Big({t_2\over t_2+t_3}\Big) 
{t_1 t_3 \over(t_1+t_3)^2}\widehat{h}({\bf t}) \ -\cr  
&- \Phi\Big({t_1 \over t_1+t_3}\Big) 
\Phi^\prime\Big({t_2\over t_2+t_3}\Big) 
{t_2 t_3 \over(t_2+t_3)^2} \widehat{h}({\bf t})  
+ \Phi\Big({t_1\over t_1+t_3}\Big) \Phi\Big({t_2\over t_2+t_3}\Big) 
{\partial [t_3\widehat{h}] \over \partial t_3} .
\end{aligned}$$
Therefore it only remains to prove that
given any point ${\bf t}^*=(t_1^*,t_2^*,0)\in \RR_+^3$
we have
\ADT
\newcounter{TtS} \setcounter{TtS}{\value{SEC}}
\newcounter{TtT} \setcounter{TtT}{\value{THM}}
$$\begin{aligned}
&D_m \chi({\bf t}) \!\!\to\! 0 \ (m \!=\!\! 1,2) ,\quad 
D_3 \chi({\bf t}) \!\!\to\! h(t_1^*,t_2^*) \cr
&\hbox{whenever} \quad  
\RR^2_+\!\times\RR_{++} \!\ni\! {\bf t} \!\to\! {\bf t}^* .
\end{aligned}
\leqno(\NUM)$$
Actually, these relation follow from the mean value expressions
$$\begin{aligned}
{\partial [t_3\widehat{h}]\over \partial t_m} &=
{1 \over t_3} \int_{s_2=t_{3-m}}^{t_{3-m}+t_3}
\big[ h(t_m+t_3,s_2) - h(t_m,s_2) \big] \, ds_2 =\cr 
&= h\big( t_m + t_3,r_{3-m}({\bf t})\big) - h\big( t_m,q_{3-m}({\bf t})\big), \cr
{\partial [t_3\widehat{h}]\over \partial t_3} &=
-{1 \over t_3^2} \int_{s_1=t_1}^{t_1+t_3} \int_{s_2=t_2}^{t_2+t_3}
h(s_1,s_2) \ d s_2 \ d s_1 \ +\cr
&\hskip5mm + {1 \over t_3} \int_{s_2=t_2}^{t_2+t_3}
h(t_1+t_3,s_2) \ d s_2 \ +
{1 \over t_3} \int_{\! s_1=t_1}^{t_1+t_3} h(s_1,t_2+t_3) \ d s_1 =\cr 
&= -h\big( p_1({\bf t}), p_2({\bf t})\big) +
h\big( t_1+t_3,r_2({\bf t})\big) +
h\big( r_1({\bf t}),t_2+t_3\big)
\end{aligned}$$
with suitable 
$$p_1({\bf t}),q_1({\bf t}), r_1({\bf t}) \in [t_1,t_1+t_3], \quad 
p_2({\bf t}),q_2({\bf t}), r_2({\bf t}) \in [t_2,t_2+t_3] .$$
Indeed, let 
$\RR_+^2\times\RR_{++} \ni {\bf t}\to {\bf t}^*=(0,t_2^*,t_3^*)$.
Then ${\bf p}({\bf t})\to (t_1^*,t_2^*)$, 
${\bf q}_m({\bf t}),{\bf r}_m({\bf t}) \to t_m^*$ $(m=1,2)$, 
whence
$\partial [t_3 \widehat{h}]/\partial t_m \to 0$ $(m=1,2)$
and
$\partial [t_3 \widehat{h}]/\partial t_3 \to h(t_1^*,t_2^*)$.
Since the function $(s,t)\!\mapsto\! s/(s \!+\! t)$ is analytic 
on the open half plane
$\{ (s,t) \!\in\!\RR^2 : s \!+\! t \!>\! 0\}$ containing the rays
$R_1 \!=\! \RR_{++} \!\times\! \{ 0\}$ 
resp. $R_2 \!=\! \{ 0\} \!\times\! \RR_{++}$, 
$(\theTtS.\theTtT)$ is immedite in the cases
${\bf t}^* \in R_m$ $(m=1,2)$.
If $\RR^2_+\times \RR_{++} \ni {\bf t} \to (0,0,0)$,
we can deduce $(\theTtS.\theTtT)$
from the facts that,
for $\ell=1,2,3$,
the functions 
$t_\ell t_3 / (t_\ell+t_3)^2, t_3^2 /(t_\ell+t_3)^2$
resp.
$\Phi\big(t_\ell / (t_\ell+t_3)\big),
\Phi^\prime\big(t_\ell / t_\ell+t_3)\big)$
are bounded,
furthermore 
$\partial[ t_3 \widehat{h}] / \partial t_\ell \to h(0,0)=0$.

%
\DEL{
$$\begin{aligned}
&\widehat{h}({\bf t}) = h\big( p_1({\bf t}),p_2({\bf t}),t_3\big) ,
\qquad p_m({\bf t}) \in [t_m,u_m] \qquad (m=1,2),\cr
&{1 \over t_3} \int_{s_2=t_2}^{t_2+t_3}
\big[ h(t_1+t_3,s_2) - h(t_1,s_2) \big] \, ds_2 =\cr 
&= {1\over t_3^2} h\big( t_1 + t_3,r_2({\bf t})\big) - 
t_3^{-1} h\big( t_1,q_2({\bf t})\big) \cr
&\int_{s_1=t_1}^{u_1} h(s_1,u_2) \, ds_1 = 
(u_1-t_1) h\big( s_1,q_1(t_1,u_1,u_2)\big) , \ \ 
q_1(t_1,u_1) \in [t_1,u_1] ,\cr  
&\int_{s_2=t_2}^{u_2} h(u_1,s_2) \, ds_2 = 
(u_2-t_2) h\big( u_1,r_2({\bf t})\big) , \ \ 
r_m({\bf t}) \in [t_m,u_m] \ \ (m=1,2)\cr
\end{aligned}$$
}


\bigskip

\ADT
\newcounter{CORiiS} \setcounter{CORiiS}{\value{SEC}}
\newcounter{CORiiT} \setcounter{CORiiT}{\value{THM}}
\noindent
{\bf \NUM. Corollary.}  {\it By Theorem $\theTHMiS.\theTHMiT$ and 
Lemma $\theLEMiiS.\theLEMiiT$, for 
any admissible pair $[\Psi_0,\Psi_1]$,
there exist $\chi_0,\chi_1 \!\in\!
{\cal C}^1_0\big(\RR_+^3\setminus\{ (0,0,0\}\big)$
such that $[\Psi_0,\Psi_1,\chi_0,\chi_1]$ is an RSD-tuple.
}

\bigskip
\ADT
\newcounter{REMvS} \setcounter{REMvS}{\value{SEC}}
\newcounter{REMvT} \setcounter{REMvT}{\value{THM}}
\noindent
{\bf \NUM. Remark.} Unfortunately,
our convolution type generic construction 
provided by the proof of Lemma $\theLEMiiS.\theLEMiiT$ 
is far from being optimal in most cases from numerical view points.

In contrast,
in view of Theorem $\theTHMiS.\theTHMiT$, in 
many cases we may apply the following satisfactory construction: 

\bigskip
\ADT
\newcounter{PROPvS} \setcounter{PROPvS}{\value{SEC}}
\newcounter{PROPvT} \setcounter{PROPvT}{\value{THM}}
\noindent
{\bf \NUM. Proposition.} {\it 
If $[\Psi_0,\Psi_1]$ is an admissible pair of shape functions
and $h_0,h_1 \!\in\!{\cal C}^1_0(\RR^2_+)$ are functions such that
$h(t,1\!-\! t) \!=\! \Psi_0^\prime(t)$ and 
$h(t,1\!-\! t) \!=\! \Psi_1^\prime(t)(1\!-\! t)$
$(0\!\le\! t\!\le\! 1)$ 
then $[\Psi_0,\Psi_1,\chi_0,\chi_1]$ is an RSD tuple with 
$$\chi_r(t_1,t_2,t_3) = h_r(t_1,t_2) t_3 \qquad
(r=0,1) .$$
}
\vskip-4mm
\ADT
\newcounter{EXxS} \setcounter{EXxS}{\value{SEC}}
\newcounter{EXxT} \setcounter{EXxT}{\value{THM}}
\noindent
{\bf \NUM. Examples.} Cases with a factorization
\vskip-7mm 
$$\Psi_0^\prime(t)= w_{01}(t)w_{02}(1-t), \quad 
\Psi_1^\prime(t)(1-t)=w_{11}(t)w_{12}(1-t)$$
\vskip-2mm \noindent
for suitable functions $w_{rk}\in {\cal C}^1_0(\RR_+)$
such that $w_{rk}(0) = w_{rk}^\prime(0) =0$. 

(a) The procedure $(\theFS.\theFT)$ 
with the 
shape functions $(\theSHS.\theSHT)$ corresponds to the case
$[\Psi_0,\Psi_1]=[\Phi,\Theta]$, 
$\chi_0(t_1,t_2,t_3)=30 t_1^2 t_2^2 t_3$,
$\chi_1(t_1,t_2,t_3)\!=\! 12 t_1^2 t_2^2 t_3$
with 
$w_{01}(t) \!=\! 30 t^2$, $w_{11}(t) \!=\! 12 t^2$, 
$w_{02}(t) \!=\! w_{12}(t) \!= \!t^2$ 
since $\Phi^\prime(t) \!=\! 30t^2(1 \!-\! t)^2$ and
$\Theta^\prime(t) \!=\! 12 t^2(1 \!-\! t)$.

\smallskip

(b) $\big[\Phi,\Phi,30t_1^2t_2^2t_3,30t_1^2t_2^3t_3\big]$ is an RSD tuple
by Proposition $\thePROPvS.\thePROPvT$,
corresponding to the case \  
$w_{01}(t)=w_{11}(t)=30 t^2$, \ $w_{02}(t)=t^2$, $w_{12}(t)=t^3$.

\bigskip

\ADT
\newcounter{REMxvS} \setcounter{REMxvS}{\value{SEC}}
\newcounter{REMxvT} \setcounter{REMxvT}{\value{THM}}
\noindent
{\bf \NUM. Remark.} It is an easy consequence of Whitney's extension theorem
[\WHIT]
that for any function $\psi\in{\cal C}([0,1])$ 
with $\psi(0)=\psi^\prime(0)=\psi^\prime(1)=0$
there exists a function
$h\in{\cal C}^1_0(\RR^2_+)$ such that $\psi^\prime(t)=h(t,0-t)$.
Hence every admissible pair $[\Psi_0,\Psi_1]$ with the property
$\Psi_0^\prime(1)=0$
admits functions $h_0,h_1\in{\cal C}(\RR^2_+)$ such that
$[\Psi_,\Psi_1,h_0(t_1,t_2)t_3,h_0(t_1,t_2)t_3]$
is an RSD tuple.

\bigskip
We complete this section with a brief description of the
local approximation operator corresponding to 
Theorem $\theTHMiS.\theTHMiT$ and its use in constructing
${\cal C}^1$-splines over a 2D triangular mesh analogously as done in
[\MDPI, Algorithm].

\bigskip
\ADT
\newcounter{DEFxS} \setcounter{DEFxS}{\value{SEC}}
\newcounter{DEFxT} \setcounter{DEFxT}{\value{THM}}
\noindent
{\bf \NUM. Definition.} Let ${\bm \Pi}= [\Psi_0,\Psi_1,\chi_0,\chi_1]$ be an
RSD tuple and write 
${\bm \Psi}=[\Psi_0,\Psi_1]$ resp. $\Chi=[\chi_0,\chi_1]$.
We introduce the {\it SDR modification of  
$\FFg^{\bm \Psi}$ by means of the complentary shape functions $\Chi$} 
as the operator
$$\FFg^{\bm \Pi} \!\! : \!  
\big\{ ({\bf p}_k,\! f_k, \! A_k,\! {\bf u}_k) \big\}_{k=1}^3 \mapsto 
\FFg^{\bm \Psi} \big\{ ({\bf p}_k,\! f_k,\! A_k) \big\}_{k=1}^3 -
\hbox{\Goth H}^\Chi \big\{ ({\bf p}_k,\! f_k,\! A_k,\!{\bf u}_k) \big\}_{k=1}^3 $$
defined for all tuples 
$\big\{ ({\bf p}_k,\! f_k,\! A_k,\!{\bf u}_k) \big\}_{k=1}^3$
where $\big\{ ({\bf p}_k,\! f_k,\! A_k) \big\}_{k=1}^3$ is a 
non-degenerate function germ and  
${\bf u}_k \not\parallel {\bf p}_i \! - {\bf p}_j$ 
$\big( (i,j,k)\in S_3\big)$ and 
$$\hbox{\Goth H}^\Chi \big\{ ({\bf p}_k,\! f_k,\! A_k,\!{\bf u}_k) \big\}_{k=1}^3 = \Big[ \hbox{ $H$ in Theorem $\theTHMiS.\theTHMiT$} \Big] .$$
Recall [\MDPI] that, given a set 
${\bf V}=\big\{ {\bf v}_1,\ldots,{\bf v}_R \big\}$ 
of points in $\RR^2$,
by a {\it triangular mesh} over ${\bf V}$
we mean a family
$\big\{ {\bf T}_1,\ldots,{\bf T}_M \big\}\subset \RR^2$ 
of non-degenerate triangles of the form
$${\bf T}_m = {\rm Conv} 
\big\{ {\bf v}_{i(m,1)},  {\bf v}_{i(m,2)}, {\bf v}_{i(m,3)} \big\}, 
\quad i(m,1)<i(m,2)<i(m,3)$$
with pairwise disjoint interiors
whose intersections are 
either empty or a common vertice
or a common edge,
furthermore ${\bf V}\subset \bigcup_{m=1}^M {\bf T}_m$.
A crutial ingreadient of the construction,
we enumerate the edges of the mesh triangles ${\bf T}_m$
in the form ${\bf E}_1,\ldots,{\bf E}_S$
and associate a vector ${\bf u}_s \not\parallel{\bf E}_s$ 
with each edge ${\bf E}_s$,
furthermore let
$s(m,k)$ denote the index of the {\it opposite} edge to 
the vertex ${\bf v}_{m(i,k)}$ in the triangle ${\bf T}_m$. 
By replacing $\Phi,\Theta$ with $\Psi_0,\Psi_1$,
respectively the functions 
$30 \lambda_i^2 \lambda_j^2 \lambda_k$, $12 \lambda_i^2 \lambda_j^2 \lambda_k$
with $\chi_0,\chi_1$ in a straightforward manner in
[\MDPI, Proof of Thm. 2], we can conclude the folowing:
Given any germ 
$\big\{ ({\bf v}_i,f_i,A_i) \big\}_{i=1}^R$ 
of first order function data over ${\bf V}$,
the union of the functions
$$F_m \!= \FFg^{\bm \Pi} \Big\{ \big({\bf v}_{i(m,k)},f_{i(m,k)},A_{i(m,k)},
{\bf u}_{s(m,k)} \big) \Big\}_{\! k=1}^{\! 3} 
\!\in {\cal C}^1({\bf T}_m) \quad (m \!=\! 1,\ldots,M)$$ 
is continuously differentiable on the mesh domain 
${\bf D} =\bigcup_{m=1}^M {\bf T}_m$.

\vskip10mm
\ADS
\noindent{\bfs \theSEC. 
 Range shift property}
\rm

\medskip
\ADT
\newcounter{DEFivS} \setcounter{DEFivS}{\value{SEC}}
\newcounter{DEFivT} \setcounter{DEFivT}{\value{THM}}
One of the most frequent applications of splines over 
triangular meshes is reconstucting approximately a smooth surface
with a function graph
passing through a family of points.
Such procedures correspond to a model of the following pattern: 
We are given a smooth surface 
${\bf S}\subset \RR^3$ which can be represented
in the form 
${\bf S}= {\rm graph}(f) = 
\big\{ \big[{\bf p},f({\bf p})\big] : {\bf p}\in\RR^2\big\}$
with some (a priori unknown) function $f\in{\cal C}^1(\RR^2)$.
Accurate data (obtained with scanner equipments) are available
for a finite family of points ${\bf v}_1,\ldots,{\bf v}_N \in\RR^2$
concerning the {\it first order data} \   
$f^{[1]}({\bf p}) = \big[ {\bf p},f({\bf p}), f^\prime({\bf p}) \big]$ \
and we approximate a piece of ${\bf S}$ with the graph of
a ${\cal C}^1$-spline function $F \in{\cal C}^1({\bf D})$ 
on the domain ${\bf D}$ of a triangular mesh with vertices ${\bf v}_n$.
The procedures used in constructing $F$
are usually investigated from the view points of estimating 
various biases between ${\rm graph}(F)$ and ${\bf S}$, 
respectively from exhibiting features of algebraic-geometric character.
In this note we only discuss two kinds of the last category in the context of 
RSD methods:
if the the spline graps are invariant with respect to tranlations of 
first order data and the property that plane surfaces remain invariant.  

\medskip
\noindent
{\bf \NUM. Definition.} 
Throughout this section  $[{\bf P,U}]$ denotes an
admissible pair of triples
${\bf P}=\big\{ {\bf p}_i \big\}_{i=1}^3,
{\bf U}=\big\{ {\bf u}_i \big\}_{i=1}^3$
with ${\bf p}_i,{\bf u}_i \in\RR^2$.
We write ${\bf T}={\rm Conv}({\bf P})$ and 
$\lambda_1,\lambda_2,\lambda_3$ resp. $G_1,G_2,G_3$ 
for the weight functions resp. their derivatives
associated with the vertices 
of the triangle ${\bf T}$.
Given any function $f\in{\cal C}^1(\RR^2)$,
we shall use the shorthand notations
\vskip-3mm 
$$f^{[1]}_{\bf P,U} = \big\{ \big[
{\bf p}_i,f({\bf p}_i),f^\prime({\bf p}_i),{\bf u}_i\big]\big\}_{i=1}^3, \quad
\hbox{\Goth F}^{\bm \Pi}_{\bf P\! ,U} f =
\hbox{\Goth F}^{\bm \Pi} f^{[1]}_{\bf P,U}  .$$
An RSD-tuple 
${\bm \Pi}=[\Psi_0,\Psi_1,\chi_0,\chi_1]$ 
has the {\it range shift property} if the procedure
$\hbox{\Goth F}^{\bm \Pi}$ leaves the constant functions invariant, 
that is if for all admissible pairs \ ${\bf P,U}$ \ 
and
the constant unit function ${\bf 1} :\RR^2 \ni {\bf x} \mapsto 1$
we have
\vskip-3mm
$$\hbox{\Goth F}^{\bm \Pi}_{\bf P,U}  {\bf 1} =
{\bf 1} \vert {\rm Conv}({\bf P})$$
\ADT
\newcounter{REMxiiiS} \setcounter{REMxiiiS}{\value{SEC}}
\newcounter{REMxiiiT} \setcounter{REMxiiiT}{\value{THM}}
\noindent
{\bf \NUM. Remark.} Since 
${\bf 1}({\bf p}_i) \!=\! 1, \  {\bf 1}^\prime({\bf p_i}){\bf u}_i \!=\! 0$,
we can write
\vskip-8mm 
$$\begin{aligned}
\hbox{\Goth F}^{{\bm \Pi}}_{\bf P \! ,U}  {\bf 1} &=
\sum_{i=1}^3 \Psi_0(\lambda_i) - \hskip-4mm \sum_{(\ell,m,n)\in{\rm S}_3}
\hskip-1mm {G_\ell {\bf u}_n \over G_n{\bf u}_n}
\chi_0(\lambda_\ell,\lambda_m,\lambda_n) =\cr 
&\hskip-8mm 
=\sum_{n=1}^3 \left\{ \Psi_0(\lambda_n) \!-\!
\left[ {G_{\ell_n} {\bf u}_n \over G_n{\bf u}_n} 
\chi_0(\lambda_{\ell_n},\lambda_{m_n},\lambda_n) \!+\!
{G_m {\bf u}_n \over G_n{\bf u}_n} 
\chi_0(\lambda_{m_n},\lambda_{\ell_n},\lambda_n) \right] \right\}
\end{aligned}$$
\vskip-2mm 
\noindent
with $(\ell_1,m_1)=(2,3)$, $(\ell_2,m_2)=(1,3)$, $(\ell_3,m_3)=(1,2)$.

\bigskip
\ADT
\newcounter{LEMixS} \setcounter{LEMixS}{\value{SEC}}
\newcounter{LEMixT} \setcounter{LEMixT}{\value{THM}}
\noindent
{\bf \NUM. Lemma.} {\it For any $\alpha_1,\alpha_2,\alpha_3\in\RR$
there exist ${\bf u}_1,{\bf u}_1,{\bf u}_1\in\RR^2$ such that
\vskip-3mm 
$$\left[ {G_i {\bf u}_j \over G_j {\bf u}_j} \right]_{i,j=1}^3 =
\left[ \begin{matrix} 
1           &-1-\alpha_2 &\alpha_3 \cr  
\alpha_1    &1           &-1-\alpha_3 \cr  
-1-\alpha_1 &\alpha_2    &1 
\end{matrix} \right] .$$
}
\vskip-3mm
{\bf Proof.} A suitable choce is  \ 
${\bf u}_1 = 
{\bf p}_1-{1\over 2}[{\bf p}_2+{\bf p}_3] + 
({1\over 2} + \alpha_1) [{\bf p}_2-{\bf p}_3]$, \ 
${\bf u}_2 = 
{\bf p}_2-{1\over 2}[{\bf p}_3+{\bf p}_1] + 
({1\over 2} + \alpha_2) [{\bf p}_3-{\bf p}_1]$, \ 
${\bf u}_3 = 
{\bf p}_3-{1\over 2}[{\bf p}_1+{\bf p}_2] + 
({1\over 2} + \alpha_3) [{\bf p}_1-{\bf p}_2]$.

\bigskip
\ADT
\newcounter{LEMivS} \setcounter{LEMivS}{\value{SEC}}
\newcounter{LEMivT} \setcounter{LEMivT}{\value{THM}}
\noindent
{\bf \NUM. Proposition.} 
{\it An RSD tuple  
${\bm \Pi}=[\Psi_0,\Psi_1,\chi_0,\chi_1]$
has range shift property 
if and only if
\vskip-6mm
\ADT
\newcounter{SYMS} \setcounter{SYMS}{\value{SEC}}
\newcounter{SYMT} \setcounter{SYMT}{\value{THM}}
$$\begin{aligned}
&\Psi_0(t)+\Psi_0(1-t) =1, \ \  \Psi_0^\prime(t)=\Psi_0^\prime(1-t) 
\hskip10mm (0 \le t \le 1) ,\cr
&\chi_0(t_1,t_2,t_3) = \chi_0(t_2,t_1,t_3) 
\quad 
\hskip23mm \big( (t_1,t_2,t_3)\in\Delta_3\big)
\end{aligned}
\leqno(\NUM)$$
\vskip-3mm
\noindent 
and \qquad
${\bf 1} =  
\sum\limits_{i=1}^3 \Psi_0(\lambda_i) + \hskip-9mm 
\sum\limits_{\hskip7mm(\ell,m,n)\in{\rm S}_3^+} 
\!\!\chi_0(\lambda_\ell,\lambda_m,\lambda_n)$.
}

\bigskip
{\bf Proof.} 
Consider the points 
${\bf x}_t = t {\bf p}_1+ (1-t){\bf p}_2$ $(0\le t\le 1)$. 
Since ${\bf 1}^\prime \equiv 0$,
${\bf 1} ({\bf x}_t)=1$,
$\lambda_1({\bf x}_t)=t$,  $\lambda_2({\bf x}_t)=1-t$,
and $\lambda_3({\bf x}_t)=0$,  
by assumption 
$$1 = {\bf 1} ({\bf x}_t) =
\hbox{\Goth F}^{\bm \Pi}_{\bf P,U} {\bf 1} \, ({\bf x}_t) =
\sum_{i=1}^3 \Psi_0\big(\lambda_i({\bf x}_t)\big) =
\Psi_0(t)+\Psi_0(1-t)+\Psi_0(0)$$
where $\Psi_0(0)=0$ due to $(\thePSIS.\thePSIT)$.
Thus in view of 
Remark $\theREMxiiiS.\theREMxiiiT$ and Lemma $\theLEMixS.\theLEMixT$,
the expression of $\hbox{\Goth F}^{{\bm \Pi}}_{\bf P \! ,U}  {\bf 1}$
has the form 
$$\hbox{\Goth F}^{{\bm \Pi}}_{\bf P \! ,U}  {\bf 1} =
\sum_{n=1}^3 \left\{ \Psi_0(\lambda_n) \!-\!
\big[ \alpha_n 
\chi_0(\lambda_{\ell_n},\lambda_{m_n},\lambda_n) \!+\!
(-1-\alpha_n) 
\chi_0(\lambda_{m_n},\lambda_{\ell_n},\lambda_n) \big] \right\} $$
where the coefficients $\alpha_1,\alpha_2,\alpha_3$ 
may assume any real value.
This is possible only if \ \
$\hbox{\Goth F}^{{\bm \Pi}}_{\bf P \! ,U}  {\bf 1} =
\sum\limits_{i=1}^3 \Psi_0(\lambda_i) +
\chi_0(\lambda_\ell,\lambda_m,\lambda_n) +
\chi_0(\lambda_\ell,\lambda_m,\lambda_n)+
\chi_0(\lambda_\ell,\lambda_m,\lambda_n)$
and
$$\chi_0(\lambda_{\ell_n},\lambda_{m_n},\lambda_n) =
\chi_0(\lambda_{m_n},\lambda_{\ell_n},\lambda_n) \qquad
(n=1,2,3) .$$ 
In particular, for $n\!=\! 3$ and a generic point 
${\bf x} \!=\!\sum\limits_{k=1}^3 t_k{\bf p}_k$ 
$(\sum_k t_k \!=\! 1,\ t_k\!\ge\! 0)$ 
of the triangle
${\rm Conv}({\bf P})$, we get
$\chi_0(t_1,t_2,t_3) \!\!=\!
\chi_0\big(\lambda_1({\bf x}),\lambda_2({\bf x}),\lambda_3({\bf x})\!\big) \!\!=\!
\chi_0\big(\lambda_2({\bf x}),\lambda_1({\bf x}),\lambda_3({\bf x})\!\big) \!\!=\!\!
\chi_0(t_2,t_1,t_3)$
which completes the proof.

\bigskip
\ADT
\newcounter{LEMviS} \setcounter{LEMviS}{\value{SEC}}
\newcounter{LEMviT} \setcounter{LEMviT}{\value{THM}}
\noindent
{\bf \NUM. Lemma.} {\it Given any smooth function $\phi$
defined on a domain $\Omega\!\subset\!\RR^3$ containing $\Delta_3$
and given any point ${\bf x}\!\in\!{\bf T}$,
for the function 
$f\!=\!\phi(\lambda_1,\lambda_2,\lambda_3)$ we have
$f^\prime({\bf x})=0$ if and only if
$$D_1 \phi({\bf t})=D_2 \phi({\bf t})=D_3 \phi({\bf t})
\quad \hbox{\rm where} \quad {\bf t} =
\big( \lambda_1({\bf x}),\lambda_2({\bf x}),\lambda_3({\bf x})\big) .$$
} 
\hskip 5mm 
{\bf Proof.} For any vector ${\bf u}\in\RR^2$ we have
$f^\prime({\bf x}) {\bf u} = 
\sum_{k=1}^3 D_k\varphi({\bf t}) G_k {\bf u}$.
Hence the statement is immediate from the relation
$\big\{ (G_1{\bf u},G_2{\bf u},G_3{\bf u}) :
{\bf v}\in\RR^2\big\} = \big\{ (\alpha_1,\alpha_2,\alpha_3) :
\sum_k \alpha_k=0 \big\}$.
$\big($Actually, given any triple $(i,j,k)\in S_3$ of indices,
by taking ${\bf u}={\bf p}_i-{\bf p}_j$ we have
$(G_i{\bf u} =1, (G_j{\bf u} =-1, (G_k{\bf u} =0$  
implying
$D_i\varphi({\bf t} - D_j\varphi({\bf t}$ if
$f^\prime({\bf x}) {\bf u}=0 \big)$. 
 
\bigskip
\ADT
\newcounter{PROPivS} \setcounter{PROPivS}{\value{SEC}}
\newcounter{PROPivT} \setcounter{PROPivT}{\value{THM}}
\noindent
{\bf \NUM. Corollary.} {\it The RSD tuple
${\bm \Pi}=
[\Psi_0,\Psi_1,\chi_0,\chi_1]$
has range shift property if and only if
$(\theSYMS.\theSYMT)$ holds and
\ADT
\newcounter{SYMsS} \setcounter{SYMsS}{\value{SEC}}
\newcounter{SYMsT} \setcounter{SYMsT}{\value{THM}}
$${\partial\Sigma\over\partial t_1}(t_1,t_2,t_3) \!=\!
{\partial\Sigma\over\partial t_2} (t_1,t_2,t_3) \!=\!
{\partial\Sigma\over\partial t_3} (t_1,t_2,t_3) \quad
\big( (t_1,t_2,t_3)\!\in\!\Delta_3\big) ,
\leqno(\NUM)$$  
$\Sigma(t_1,t_2,t_3) \!=\! 
\sum\limits_{i=1}^3 \!\Psi_0(t_i) +\hskip-5mm 
\sum\limits_{(\ell,m,n)\in S_3^+} \hskip-5mm\chi_0(t_\ell,t_m,t_n)$, \ 
$S_3^+ \!=\! \{ (1,\!2,\!3),(2,\!3,\!1),(3,\!1,\!2)\}$.
}

\bigskip
{\bf Proof.} 
We established that
${\bm \Pi}$ has range shift property if and only if
\ADT
\newcounter{SYMtS} \setcounter{SYMtS}{\value{SEC}}
\newcounter{SYMtT} \setcounter{SYMtT}{\value{THM}}
$$\Sigma\big(\lambda_1({\bf x}),\lambda_2({\bf x}),\lambda_3({\bf x})\big) =1 
\qquad \big({\bf x} \!\in\!{\rm Conv}({\bf P})\big) .
\leqno(\NUM)$$
Since ${\bm \Pi}$ is an RSD tuple, in any case
we have $\hbox{\Goth F}^{\bm \Pi}_{\bf P,U} ({\bf p}_i) =1$ $(i=1,2,3)$.
Hence the relation $(\theSYMtS.\theSYMtT)$ holds if and only if
the Fr\'echet derivative of the function 
$\hbox{\Goth F}^{\bm \Pi}_{\bf P,U} {\bf 1}$ vanishes, that is if 
$$\big[\Sigma(\lambda_1,\lambda_2,\lambda_3)\big]^\prime({\bf x}) =0 
\qquad \big({\bf x} \!\in\!{\rm Conv}({\bf P})\big) . 
\leqno(\theSYMtS.\theSYMtT^\prime)$$
An application of Lemma $\theLEMviS.\theLEMviT$ with
$\phi=\Sigma$ completes the proof. 

\DEL{
Notice that, in terms of the weight values $t_k=\lambda_k({\bf x})$
we can write
$$\begin{aligned}
&\big[\Sigma(\lambda_1,\lambda_2,\lambda_3)\big]^\prime({\bf x}) =
\sum_{k=1}^3 
{\partial\Sigma \over \partial t_k}(t_1,t_2,t_3) G_k =\cr 
&= {\partial\Sigma \over \partial t_1}(t_1,t_2,t_3) G_1 +
{\partial\Sigma \over \partial t_2}(t_1,t_2,t_3) G_2 +
{\partial\Sigma \over \partial t_3}(t_1,t_2,t_3) \big[-G_1-G_2\big] =\cr
&= \Big[ {\partial\Sigma \over \partial t_1}(t_1,t_2,t_3) \!-\!
{\partial\Sigma \over \partial t_3}(t_1,t_2,t_3) \Big] G_1 \!+\!
\Big[ {\partial\Sigma \over \partial t_2}(t_1,t_2,t_3) \!-\!
{\partial\Sigma \over \partial t_3}(t_1,t_2,t_3) \Big] G_2 .
\end{aligned}$$
Since $\big\{ \big( \lambda_1({\bf x}),\lambda_2({\bf x}),\lambda_3({\bf x}) :
{\bf x} \in {\rm Conv}({\bf P}) \big\} = \Delta_3$
and since any pair
$G_i=\lambda_i^\prime$, $G_j=\lambda_j^\prime$ $(1\le i< j\le 3)$ of
derivative weights is linearly independent 
(but $G_1+G_2+G_3=0$),  
we can conclude that relations $(\theSYMsS.\theSYMsT),(\theSYMtS.\theSYMtT)$
and $(\theSYMtS.\theSYMtT^\prime)$ are equivalent.
}

\bigskip
\ADT
\newcounter{EXivS} \setcounter{EXivS}{\value{SEC}}
\newcounter{EXivT} \setcounter{EXivT}{\value{THM}}
\noindent
{\bf \NUM. Example.} The RSD procedures corresponding to tuples of the form
$\big[ t^3(10-15t+6t^2),*,30t_1^2t_2^2t_3,*\big]$
as in Examples $\theEXxS.\theEXxT$ have range shift property.
This highly non-trivial fact can easily be established by verifying
$(\theSYMsS.\theSYMsT)$ as follows.
In such cases, we have

\centerline{
$\Sigma(t_1,t_2,t_3) \!=\! \sum\limits_{i=1}^3 t_i^3(10-15t_i+6t_i^2) + 
30 \big( t_1^2 t_2^2 t_3 + t_2^2 t_3^2 t_1 + t_3^2 t_1^2 t_2)$
}

\noindent
and it suffices to show that the polynomial
${\partial\Sigma \over \partial t_1}$ 
is symmetric
in $t_1,t_2,t_3$ when replacing $t_1=1-t_2-t_3$.
Since $\big[ t^3 (10-15t+6t^2)\big]^\prime = 30 t^2(1-t)^2$,
if $(t_1,t_2,t_3)\in\Delta_3$ then 
$\big[ t_1^3 (10-15t_1+6t_1^2)\big]^\prime = 30 t_1^2(t_2+t_3)^2$ and 
$$\begin{aligned}
&{1\over 30}{\partial\Sigma \over \partial t_1} = 
t_1^2(t_2+t_3)^2 + 2 t_1 t_2^2 t_3 + t_2^2 t_3^2 + 2 t_3^2 t_1 t_2 =\cr 
&= t_1^2 t_2^2 + t_1^2 t_3^2 + t_2^2 t_3^2 + 2 t_1 t_2 t_3 (t_1+t_2+t_3) =
 t_1^2 t_2^2 + t_1^2 t_3^2 + t_2^2 t_3^2 + 2 t_1 t_2 t_3 .
\end{aligned}$$

\ADT
\newcounter{DEFvS} \setcounter{DEFvS}{\value{SEC}}
\newcounter{DEFvT} \setcounter{DEFvT}{\value{THM}}
\bigskip
\noindent
{\bf \NUM. Definition.} We say that a function
$\Psi\in{\cal C}^1([0,1])$ is {\it d-symmetric} if
\vskip-3mm
$$\Psi^\prime(t) = \Psi^\prime(1-t) \quad (0\le t\le 1) .$$

\ADT
\newcounter{CORiiiS} \setcounter{CORiiiS}{\value{SEC}}
\newcounter{CORiiiT} \setcounter{CORiiiT}{\value{THM}}
\noindent
{\bf \NUM. Corollary.} {\it 
If ${\bm \Pi}=[\Psi_0,\Psi_1,\chi_0,\chi_1]$ is an RSD tuple
with d-symmetric shape function $\Psi_0$
then 
${\bm \Pi}^{[s]}=[\Psi_0,\Psi_1,\chi_0^{[s]},\chi_1]$
with
\ADT
\newcounter{CHInsS} \setcounter{CHInsS}{\value{SEC}}
\newcounter{CHInsT} \setcounter{CHInsT}{\value{THM}} 
$$\chi_0^{[s]}(t_1,t_2,t_3) = 
{1\over 2} \chi_0(t_1,t_2,t_3) + {1\over 2} \chi_0(t_2,t_1,t_3)
\leqno(\NUM)$$
is also an RSD tuple such that
\ADT
\newcounter{FFS} \setcounter{FFS}{\value{SEC}}
\newcounter{FFT} \setcounter{FFT}{\value{THM}}
$$\hbox{\Goth F}^{{\bm \Pi}^{[s]}}_{\bf P,U} {\bf 1} =
\sum_{(\ell,m,n)\in S_3} \Big\{ {1\over 2} \Psi_0(\lambda_\ell) + 
\chi^{[s]}(\lambda_\ell,\lambda_m,\lambda_n)\Big\} .
\leqno(\NUM)$$
}

{\bf Proof.}
The marginal conditions $(\theCHImS.\theCHImT)$ for $\chi_0^{[s]}$
are immediate.  
As a consequence of the d-symmetry of $\Psi_0$, we have
$$\begin{aligned}
D_3\chi_0^{[s]}(t,1-t,0) &= 
{1\over 2} D_3\chi_0(t,1-t,0) + {1\over 2} D_3\chi_0(1-t,t,0) =\cr
&= {1\over 2} \Psi_0(t) + {1\over 2} \Psi_0(1-t) = \Psi_0(t)
\end{aligned}$$
establishing that ${\bm \Pi}^{[s]}$ is an RSD tuple.
Due to the symmetry of $\chi_0^{[s]}$ in the variables $t_1,t_2$
in Remark $\theREMxiiiS.\theREMxiiiT$ applied to 
$\chi_0^{[s]}$ we can write
$$\hbox{\Goth F}^{{\bm \Pi}}_{\bf P \! ,U}  {\bf 1} =
\sum_{n=1}^3 \left\{ \Psi_0(\lambda_n) \!-\!
\left[ \Big[ {G_{\ell_n} {\bf u}_n \over G_n{\bf u}_n} 
\!+\! {G_m {\bf u}_n \over G_n{\bf u}_n} \Big] 
\chi_0(\lambda_{\ell_n},\lambda_{m_n},\lambda_n) \right] \right\} .$$ 
Taking Lemma $\theLEMixS.\theLEMixT$ into account,
we obtain $(\theFFS.\theFFT)$.

\DEL{
\ADT
\newcounter{LEMvS} \setcounter{LEMvS}{\value{SEC}}
\newcounter{LEMvT} \setcounter{LEMvT}{\value{THM}}
\bigskip
\noindent
{\bf \NUM. Lemma.} {\it 
Assume $h\in{\cal C}_0^1({\bf T})$ 
with
$h^\prime({\bf x})=0$ 
for ${\bf x}\in\partial{\bf T}$.
Then the function 
$\widehat{h}$ defined for 
${\bf t}=(t_1,t_2,t_3) \in \RR^3_+$ by
\ADT
\newcounter{HAThS} \setcounter{HAThS}{\value{SEC}}
\newcounter{HAThT} \setcounter{HAThT}{\value{THM}} 
$$\begin{aligned}
&\widehat{h}({\bf t}) :=
{\rm s}({\bf t}) h\big( {\rm s}({\bf t})^{-1} {rm S}({\bf t}) \big) 
\quad (t_3>0), \quad \widehat{h}(0,0,0):=0 \cr  
&\quad \hbox{\rm where} \ \ 
{\rm s}({\bf t})=t_1+t_2+t_3, \ \ 
{\rm S}({\bf t})=t_1{\bf p}_1+t_2{\bf p}_2+t_3{\bf p}_3
\end{aligned}
\leqno(\NUM)$$
belongs to ${\cal C}_0^1([0,1]^3)$ with the marginal behaviour
$$\widehat{h} (t_1,t_2,t_3) = 
D_m \widehat{h} (t_1,t_2,t_3) =0 \quad \hbox{\rm for} \ \ 
m=1,2,3 \ \ \hbox{\rm if} \ \ t_1 t_2 t_3 =0$$
and being such that in terms of the weights of ${\bf T}$ we have 
$$\widehat{h}
\big(\lambda_1({\bf x}),\lambda_2({\bf x}),\lambda_3({\bf x})\big) =
h({\bf x}) \qquad ({\bf x}\!\in\!{\bf T}) .$$
}
{\bf Proof.} We have to see only that,
in terms of Fr\'echet derivatives, 
$\widehat{h}^\prime({\bf t})=0$ for 
${\bf t}\in\partial\RR^3_+ =
\big\{ (t_1,t_3,t_3)\in\RR^3_+ : t_1=0 \ \hbox{or} \
t_2=0 \ \hbox{or} \ t_3=0 \big\} = 
\big\{ (t_1,t_2,t_3): t_1,t_2,t_3\ge 0, \ t_1 t_2 t_3 =0\big\}$.
Notice that 
$\partial{\bf T} = 
\big\{ {\rm S}({\bf t}) : {\bf t} \in \partial\RR_+^3\big\}$.
Thus given 
${\bf t}=(t_1,t_2,t_3) \in \partial\RR^3_+$
with $t_3>0$,
the point 
${\bf x}: = {\rm s}({\bf t})^{-1} {\rm S}({\bf t})$
belongs to 
$\partial{\bf T}$ 
and therefore
$$\widehat{h}^\prime({\bf t}) = 
\big[ {\rm s} \cdot h\big( {\rm s}^{-1} {\rm S} \big) \big]^\prime({\bf t}) =
h({\bf x}) \cdot \big[- {\rm s}({\bf t})^{-2} {\rm s}^\prime({\bf t})\big] + 
{\rm s}({\bf t})^{-1} h^\prime ({\bf x}) =0 $$
which completes the proof.
}

\bigskip
\ADT
\newcounter{THMiiS} \setcounter{THMiiS}{\value{SEC}}
\newcounter{THMiiT} \setcounter{THMiiT}{\value{THM}}
\noindent
{\bf \NUM. Theorem.} {\it
Let $[\Psi_0,\Psi_1]$ be an admissible pair of shape functions,
with $\Psi_0$ being d-symmetric. 
Then we can find $\chi_0^*,\chi_1^*\in{\cal C}_0(\RR^3_+)$ 
such that ${\bm \Pi}^*=[\Psi_0,\Psi_1,\chi_0^*,\chi_1^*]$ is an RSD-tuple
with range shift property 
where
$\chi_0^*$ is symmetric 
in its first two variables.
}

\bigskip
 
{\bf Proof.} According to Theorem $\theTHMiS.\theTHMiT$
there exist an RSD tuple of the form
${\bm \Pi}=[\Psi_0,\Psi_1,\chi_0,\chi_1]$
with $\chi_0,\chi_1\in{\cal C}_0^1(\RR_0^3)$. 
By Corollary $(\theCORiiiS.\theCORiiiT)$, 
${\bf \Pi}^{[s]} =[\Psi_0,\Psi_1,\chi_0^{[s]},\chi_1]$
is also an RSD tuple 
where the complementary shape function 
$\chi_0^{[s]}$ is symmetric in its first two variables
and such that $(\theFFS.\theFFT)$ holds.
Then, independently of the choice of the vector triple
${\bf U}=[{\bf u}_1,{\bf u}_2,{\bf u}_3]$, 
we can write 
the difference function of 
${\bf 1}$ and
$\hbox{\Goth F}^{{\bm \Pi}^{[s]}}_{\bf P \! ,U} {\bf 1}$
in the form 
\ADT
\newcounter{FFdS} \setcounter{FFdS}{\value{SEC}}
\newcounter{FFdT} \setcounter{FFdT}{\value{THM}}
$$\begin{aligned}
&d = {\bf 1} -
\hbox{\Goth F}^{{\bm \Pi}^{[s]}}_{\bf P \! ,U} {\bf 1} = 
\delta(\lambda_1,\lambda_2,\lambda_3)
\quad \hbox{with} \cr 
&\delta({\bf t}) = 
d\big( \hbox{$\sum_{k=1}^3 t_k{\bf p}_k \big)$} \!=
\hskip-4mm \sum_{(\ell,m,n)\in{\rm S}_3} \!\!\! 
\left\{ {1\over 6} -
{1\over 2} \Psi_0\big(t_\ell\big) - 
\chi_0^{[s]}\big( t_\ell,t_m,t_n\big) \right\} .
\end{aligned}
\leqno(\NUM)$$
Observe that 
\ADT
\newcounter{DeltaS} \setcounter{DeltaS}{\value{SEC}}
\newcounter{DeltaT} \setcounter{DeltaT}{\value{THM}}
$$d({\bf x})=0, \quad 
d^\prime({\bf x})=0 \qquad ({\bf x}\in\partial{\bf T}) .
\leqno(\NUM)$$
Indeed, at the generic point $(\theGENPS.\theGENPT)$ of the edge opposite to
the vertex ${\bf p}_k$,
by $(\theFOS.\theFOT)$ we have 
$d({\bf x}^k_t) \!=\! 1 \!-\! [\Psi_0(t) \!+\! \Psi_0(1\!-\! t)] \!=\! 0$.
In particular the functions
${\bf 1}$ and $\hbox{\Goth F}^{\Pi^{[s]}}_{\bf P,U}$
coincide along the edge $[{\bf p}_i,{\bf p}_j]$,
thus
$d^\prime({\bf x}_t^k) ({\bf p}_j \!-\! {\bf p}_i) \!=\! 0$ with 
$d({\bf x}) \!=\! 0$ for ${\bf x}\!\in\!\partial{\bf T}$.
On the other hand, since trivially ${\bf 1}^\prime({\bf p}) \!=\!0$
everywhere,
the RSD construction ensures that 
$d^\prime({\bf x}^k_t){\bf u}_k \!=\! 
{\bf 1}^\prime({\bf x}_t^k){\bf u}_k \!-\! 
\big[\Psi_0(t) 1^\prime({\bf p}_k)){\bf u}_k \!+\!
\Psi_0(1 \!-\! t) 1^\prime({\bf p}_k)){\bf u}_k \big] \!=\! 0$.
Thus $\delta^\prime({\bf x}_t^k)$ vanishes in two
linearly independent directions implying that
$d^\prime({\bf x}_t^k){\bf u} =0$ $({\bf u}\in\RR^2)$
i.e. $d^\prime({\bf x}^k_t) =0$.

Notice that ${\bf x}\in\partial{\bf T}$ if and only if
$\big( \lambda_1({\bf x}),\lambda_2({\bf x}),\lambda_3({\bf x}) \big) 
\in \bigcup_{k=1}^3 \Delta_{3,k}$.
Hence, due to $(\theDeltaS.\theDeltaT)$,
we can apply Lemma $\theLEMviS.\theLEMviT$ with
$f=d$ and $\phi=\delta$ to conclude that
$$\delta({\bf t})=0, \quad 
D_1 \delta({\bf t}) = D_2 \delta({\bf t}) = D_3 \delta({\bf t}) \quad
\hbox{whenever} \ \ {\bf t} \in 
\hbox{$\bigcup\limits_{k=1}^3 \Delta_{3,k}$} .
\leqno(\theDeltaS.\theDeltaT^\prime)$$
On the domain $\Omega:=\RR_+^3\setminus\{ 0\}$, define
$$\begin{aligned}
&\chi_0^* = \chi_0^{[s]} + {1\over 6}\widehat{\delta} 
\quad  \hbox{where} \cr
&\widehat{\delta}({\bf t}) = \delta\Big(
{t_1\over t_1\!+\! t_2\!+\! t_3},{t_2\over t_1\!+\! t_2\!+\! t_3},
{t_3\over t_1\!+\! t_2\!+\! t_3}\Big) .
\end{aligned}$$
By writing \ $s({\bf t}) = t_1+t_2+t_3$ for short, 
for  $k=1,2,3$ we have
$$\begin{aligned}
D_k &\widehat{\delta}({\bf t}) =
{\partial \over \partial t_k} \delta\Big( {t_1 \over s({\bf t})},
{t_2 \over s({\bf t})},{t_3 \over s({\bf t})} \Big) =\cr
&= D_1\delta\Big( {{\bf t}\over s({\bf t})}\Big) 
{\partial \over \partial t_k} \, {t_1 \over s({\bf t})} +\! 
D_2\delta\Big( {{\bf t}\over s({\bf t})}\Big) 
{\partial \over \partial t_k} \, {t_2 \over s({\bf t})} +\!
D_3\delta\Big( {{\bf t}\over s({\bf t})}\Big) 
{\partial \over \partial t_k} \, {t_3 \over s({\bf t})} . 
\end{aligned}$$
In particular, 
if 
$0\ne {\bf t} \in \RR^3_+$ 
then $s({\bf t})^{-1}{\bf t} \in \Delta_3$,
furthemore
$s({\bf t})^{-1}{\bf t} \in \bigcup_{k=1}^3 \Delta_{3,k}$
whenever 
$0\ne {\bf t} \in \partial\RR^3_+$.  
Thus, as a consequence of $(\theDeltaS.\theDeltaT)$,
we get
$$
D_k \widehat{\delta} ( {\bf t} ) 
\!=\!
\sum_{\ell=1}^3 D_1 \delta \Big( {{\bf t}\over s({\bf t})} \Big)
{\partial\over \partial t_k} {t_\ell \over s({\bf t})} 
\!=\! D_1 \delta \Big( {{\bf t}\over s({\bf t})} \Big)
{\partial\over \partial t_k} {s({\bf t}) \over s({\bf t})} \!=\! 0 
\quad \big( 0\!\ne\! {\bf t} \!\in\! \partial\RR^3_+ \big)$$
It follows $D_k \chi_0^* ({\bf t}) = D_k D_k \chi_0 ({\bf t})$
for $0 \ne {\bf t} \in \partial\RR^3_0 = \big\{ (t_1,t_2,t_3)\in\RR_+^3 :
t_1 t_2 t_3 =0 \big\}$ and $k=1,2,3$
whence, in view of Theorem $\theTHMiS.\theTHMiT$ we
concluce that ${\bm \Pi}^*$ with $\chi_1^* = \chi_1$ is
indeed an RSD tuple.

\DEL{
Since $d(\cdot)$ and $d^\prime(\cdot)$
vanish along the edges of the triangle ${\bf T}$, 
according to Lemma $\theLEMvS.\theLEMvT$ we have
$\chi_0 \in {\cal C}_0(\RR^3_+)$ and 
$D_m \chi_0 = D_m \chi_0^{[s]}$ $(m=1,2,3)$.
That is $\chi_0$ has also the properties
$(\theCHIS.\theCHIT),(\theCHImS.\theCHImT)$, thus
${\bm \Pi}=[\Psi_0,\Psi_1,\chi_0,\chi_1]$ is an RSD tuple as well.
By Lemma $\theLEMvS.\theLEMvT$, also 
$$\widehat{\delta}({\bf x}) = \delta({\bf x}) 
\quad \big({\bf x} \in {\bf T}\big) .$$ 
Hence, in view of $(\theFFS.\theFFT)$ and $(\theFFdS.\theFFdT)$, 
for any point ${\bf x}\in {\bf T}={\rm Conv}({\bf P})$ 
we get
$$\hbox{\Goth F}^{\bm \Pi}_{\bf P \! ,U} {\bf 1} ({\bf x}) =
\hbox{\Goth F}^{{\bm \Pi}^{[s]}}_{\bf P \! ,U} {\bf 1} ({\bf x}) + 
\delta({\bf x}) = 
\hbox{\Goth F}^{{\bm \Pi}^{[s]}}_{\bf P \! ,U} {\bf 1} ({\bf x}) +1 - 
\hbox{\Goth F}^{{\bm \Pi}^{[s]}}_{\bf P \! ,U} {\bf 1} ({\bf x}) = 1$$
which completes the proof.

}

\DEL{ 
\ADT
\newcounter{REMxiS} \setcounter{REMxiS}{\value{SEC}}
\newcounter{REMxiT} \setcounter{REMxiT}{\value{THM}}
\bigskip
\noindent
{\bf \NUM. Remark.} In the setting of Proposition $\thePROPvS.\thePROPvT$
and in the presence of the {\it symmetry in $t_1,t_2$ of $\chi_0$},
we necessarily have 
$$\chi_0(t_1,t_2,t_3)=w_{01}(t_1)w_{02}(t_2)t_3 \ \ 
\hbox{where} \ \ 
w_{01} = \kappa w_{02}$$ 
with some constant $0\ne\kappa\in\RR$,
furthermore
$$\big[\Psi_0(\lambda_k)\big]^\prime = 
\Psi_0^\prime(\lambda_k) G_k = w_{01}(\lambda_k)w_{02}(1-\lambda_k) G_k .$$ 
Using the the abbreviation
$w = w_{01}$,
the function
$F \!= F_{\bf P,\! U} ={\hbox{\Goth F}}^{\bm \Pi}_{\bf P,U} {\bf 1}$
can be written in the form
$$
\sum\limits_{k=1}^3 \! \Psi_0(\lambda_k) -
\hskip-3mm \sum\limits_{\hskip-5mm(\ell,m,n)\in S_3} \hskip-2.2mm 
{G_\ell{\bf u}_n\over G_n{\bf u}_n}
\kappa w(\lambda_\ell)w(\lambda_m)\lambda_n 
\!=\! \sum_{k=1}^3 \! \Big[ \Psi_0(\lambda_k) + \hskip-4mm 
\sum_{i,j\!\ne\! k,\ i\!<\! j} \hskip-5mm
\kappa w(\lambda_1)w(\lambda_j)\lambda_k\Big]
$$ 
as a consequence of the relations $(\theGUS.\theGUT)$.

Since here we automatically have $F({\bf p}_i)=1$ $(i=1,2,3)$,
the procedure
${\hbox{\Goth F}}^{\bm \Pi}$ has range shift property if and only if
$F^\prime=0$ that is if, in any case, 
$F^\prime({\bf x})[{\bf p}_i-{\bf p}_j]=0$.
Assuming $i \!=\! 1, j\!=\! 2$ without loss of generality,
since $G_k[{\bf p}_1-{\bf p}_2] =\delta_{k1}-\delta_{k2}$, \ \
{\it
${\hbox{\Goth F}}^{\bm \Pi}$ has range shift property if and only if}
\ADT
\newcounter{WS} \setcounter{WS}{\value{SEC}}
\newcounter{WT} \setcounter{WT}{\value{THM}}
$$\begin{aligned}
&(\NUM) \hskip 10mm
0 =\kappa  F^\prime({\bf x})[{\bf p}_1-{\bf p}_2] =\cr
&= w(\lambda_1)w(1 \!-\! \lambda_1) \!-\! w(\lambda_2)w(1 \!-\!\lambda_2) 
\!+\! \Big[ w^\prime(\lambda_1)w(\lambda_2)\lambda_3 \!-\!
w(\lambda_1)w^\prime(\lambda_2)\lambda_3 \Big] +\cr
&+ \Big[ w(\lambda_2)w(\lambda_3) -
w^\prime(\lambda_2)w(\lambda_3)\lambda_1 \Big] +
\Big[ w^\prime(\lambda_1)w(\lambda_3)\lambda_2 -
w(\lambda_1)w(\lambda_3) \Big] . 
\end{aligned} $$
In particular, at the points 
${\bf s}_t=t{\bf p}_1+t{\bf p}_2+(1-2t){\bf p}_3$
with $\lambda_1({\bf s}_t)=\lambda_2({\bf s}_t)=t$ and
$\lambda_3({\bf s}_t)=1-2t$ 
we simply get
$$0 = F^\prime({\bf s}_t)[{\bf p}_1-{\bf p}_2] =
w(t) - t w^\prime(t), \quad w(t)= {\rm cost}\cdot t^2 \qquad (0<t<1) .$$
The relations $\Psi_0^\prime(t)=\kappa w(t)w(1-t)$, 
$\Psi_0(t)=0$, $\Psi_1(t)=1$  
imply that necessarily $\kappa=30$, $w(t)=t^2$ and 
$\chi_0(t_1,t_2,t_3) 30 t_1^2 t_2^2 t_3$. 
One can verify with straightforward calculation
(a {\tt simplify} command in computer algebra suffices) that
$(\theWS.\theWT)$ holds for $w(t)=t^2$. Thus we obtained
the conclusion:

\bigskip
\ADT
\newcounter{PROPviS} \setcounter{PROPviS}{\value{SEC}}
\newcounter{PROPviT} \setcounter{PROPviT}{\value{THM}}
\noindent{\bf \NUM. Proposition.}
{\it An RSD tuple ${\bm \Pi} \!=\! [\Psi_0,\Psi_1,\chi_0,\chi_1]$
with 
$\chi_0(t_1,t_2,t_3) \!=\!\chi_0(t_2,t_1,t_3) \!=\! w_{01}(t_1)w_{02}(t_2)t_3$ 
has range shift property if and only if
$$\Psi_0(t)=\int_0^t \!\! 30 s^2(1-s)^2 ds  \!=\!
t^3(10 \!-\! 15t \!+\! 6t^2) \ \hbox{\rm and} \
\chi_0(t_1,t_2,t_3) \!=\! 30 t_1^2 t_2^2t_3 .$$ 
In particular both procedures in {\rm Examples $\theEXxS.\theEXxT$}
have range shift property. 
}

}


\vskip10mm
\ADS
\noindent{\bfs \theSEC. 
 Affinity invariance}
\rm

\bigskip
\ADT
\newcounter{DEFviS} \setcounter{DEFviS}{\value{SEC}}
\newcounter{DEFviT} \setcounter{DEFviT}{\value{THM}}
\noindent
{\bf \NUM. Definition.} An RSD tuple 
${\bm \Pi} =[\Psi_0,\Psi_1,\chi_0,\chi_1]$
(and the procedure $\hbox{\Goth F}^{\bm \Pi}$) 
is
{\it affinity invariant} if 
$$\hbox{\Goth F}^{\bm \Pi}_{\bf P \! ,U} f = f\vert {\rm Conv}({\bf P})$$ 
for all admissible pairs $[{\bf P,\! U}]$
$\big($i.e. ${\rm Conv}({\bf P})$ is a non-degenerate triangle,
\break
${\bf p}_\ell \!-\! {\bf p}_m \!\not\parallel\! {\bf p}_n$ 
$\big( (\ell,m,n)\in S_3\big)\big)$ 
and 
affine (linear+costant) functions $f \!\! : \!\RR^2\!\!\to\!\RR$.

\bigskip
\ADT
\newcounter{REMviS} \setcounter{REMviS}{\value{SEC}}
\newcounter{REMviT} \setcounter{REMviT}{\value{THM}}
\noindent
{\bf \NUM. Remark.} (a) Affinity invariance implies range shift property.

\smallskip
(b) Given any non-degenerate triangle 
${\bf T} \!=\! {\rm Conv}\{ {\bf p}_1,{\bf p}_2,{\bf p}_3\} \!\subset\!\RR^2$,
any affine function $\RR^2\to\RR$ is a linear combination of 
the weights $\lambda_1,\lambda_2,\lambda_3$.
Hence to verify the affinity invariance of ${\bf \Pi}$, 
it suffices to prove the relation \ \ 
$\hbox{\Goth F}^{\bm \Pi}_{\bf P \! , U} \lambda_1 ({\b x}) 
=\lambda_1 ({\bf x}) \ \ \big({\bf x}\in{\rm Conv}({\bf P})\big)$ \ \
for all admissible pairs $[{\bf P,\! U}]$.

\smallskip
(c) In view of Remark $\theREMiS.\theREMiT$(c) and the 
identities $G_1({\bf x}-{\bf p}_k) = \lambda({\bf x})-\delta_{1k}$ resp.
$G_i({\bf p}_j-{\bf p}_k) = \lambda_i({\bf p}_j) - \lambda_i({\bf p}_k) =
\delta_{ij}-\delta_{ik}$,
we can write
$$\begin{aligned}
&{\hbox{\Goth F}}^{\bf \Pi}_{\bf P,\! U} \lambda_1 =
\sum_{k=1}^3 \Big[\lambda_1({\bf p}_k) \Psi_0(\lambda_k) +
G_1 ({\bf x}-{\bf p}_k) \Psi_1(\lambda_k) \Big] -\cr 
&\hskip3mm - \hskip-3mm \sum_{(\ell,m,n)\in S_3} 
\hskip-1mm {G_\ell{\bf u}_n \over G_n{\bf u}_n}
\Big[ \lambda_1({\bf p}_\ell) \chi_0(\lambda_\ell,\lambda_m,\lambda_n) +
G_1({\bf p}_m \!-\! {\bf p}_\ell) \chi_1(\lambda_\ell,\lambda_m,\lambda_n) \Big] =\cr 
&\hskip5mm = \Psi_0(\lambda_1) + \Psi_1(\lambda_1)(\lambda_1 - 1) 
+ \Psi_1(\lambda_2)\lambda_1 + \Psi(\lambda_3)\lambda_1 -\cr 
&\hskip7mm-\! {G_1 {\bf u}_3 \over G_3 {\bf u}_3} \Big\{ 
\chi_0(\lambda_1,\lambda_2,\lambda_3) \!-\! 
\chi_1(\lambda_1,\lambda_2,\lambda_3)\Big\} \!-\! 
{G_2 {\bf u}_3 \over G_3 {\bf u}_3}  
\chi_1(\lambda_2,\lambda_1,\lambda_3) - \cr
&\hskip7mm-\! {G_1 {\bf u}_2 \over G_2 {\bf u}_2} \Big\{ 
\chi_0(\lambda_1,\lambda_3,\lambda_2) \!-\! 
\chi_1(\lambda_1,\lambda_3,\lambda_2)\Big\} \!-\!
{G_3 {\bf u}_2 \over G_2 {\bf u}_2}
\chi_1(\lambda_3,\lambda_1,\lambda_2) \!-\! 0 \!-\! 0 .
\end{aligned}$$
with the ordering 
$(1,\! 2,\! 3),(2,\! 1,\! 3),(1,\! 3,\! 2),
(3,\! 1,\! 2),(2,\! 3,\! 1),(3,\! 2,\! 1)$
for  ${\rm S}_3$.

\bigskip
\ADT
\newcounter{PROPiS} \setcounter{PROPiS}{\value{SEC}}
\newcounter{PROPiT} \setcounter{PROPiT}{\value{THM}}
\noindent
{\bf \NUM. Lemma.} {\it 
If the RSD tuple  
${\bm \Pi}=[\Psi_0,\Psi_1,\chi_0,\chi_1]$
is affinity invariant then necessarily
\ADT
\newcounter{AFS} \setcounter{AFS}{\value{SEC}}
\newcounter{AFT} \setcounter{AFT}{\value{THM}}
$$\begin{aligned}
&\Psi_1(t)=\Psi_0(t) 
 \ \  \hbox{\rm and} \ \ \Psi_0(t)+\Psi_0(1-t)=1 
\quad (0\le t\le 1) ,\cr  
&\chi_0(t_1,t_2,t_3) = 
\chi_1(t_1,t_2,t_3) + \chi_1(t_2,t_1,t_3) \hskip7mm
\big( (t_1,t_2,t_3)\in \Delta_3\big) .
\end{aligned}
\leqno(\NUM)$$ 
}  
{\bf Proof.}
Since ${\bm \Pi}$ is an RSD tuple, 
by Definition $\theDEFiiS.\theDEFiiT$,
the directional derivatives
$F^\prime({\bf x}^3_t){\bf u}_3$ along the side $[{\bf p}_2,{\bf p}_3]$
with generic point ${\bf x}^1_t = t{\bf p}_2 + (1-t){\bf p}_3$ 
of the function
$F=\hbox{\Goth F}^{\bm \Pi}_{\bf P,\! U} \lambda_1$ 
are the linear combinations 
$$\begin{aligned}
F^\prime({\bf x}_t^1) {\bf u}_1 &=
\Psi_1(t)  F^\prime({\bf p}_2) {\bf u}_1  + 
\Psi_1(1-t)  F^\prime({\bf p}_3) {\bf u}_1  =\cr
&= \Psi_1(t)  G_1 {\bf u}_1  + \Psi_1(1-t)  G_1 {\bf u}_1  = 
\big[ \Psi_1(t) + \Psi_1(1-t) \big]  G_1 {\bf u}_1 .
\end{aligned}$$
Due to the admissibility of $[{\bf P,\! U}]$, $G_1 {\bf u}_1 \ne 0$.
In case of the affinity invariance of ${\bm \Pi}$,
$F=\lambda_1$ 
with $F^\prime({\bf x}_t^1) {\bf u}_1 =G_1{\bf u}_1$
implying
$\Psi_1(t) + \Psi_1(1-t) =1$. 
Also if
$F =\hbox{\Goth F}^{\bm \Pi} \lambda_1=\lambda_1$ then,
as a consequence of the relations $(\thePARS.\thePART)$,
at the points 
${\bf x}^3_t = t{\bf p}_1 + (1-t){\bf p}_2$
we have   
$$\begin{aligned}
t &= \lambda_1({\bf x}^3_t) = 
F({\bf x}^3_1) =
\Psi_0(t) - \Psi_1(t) (1-t) + \Psi_1(1-t)t, \cr 
&= \big[\Psi_0(t)-\Psi_1(t) \big] + t \big[ \Psi_1(t)+\Psi_1(1-t) \big] =
\big[\Psi_0(t)-\Psi_1(t) \big] + t
\end{aligned}$$ 
whence the first statement in $(\theAFS.\theAFT)$ is immediate.

The proof of the second statement relies upon the fact that
the formula for $\hbox{\Goth F}^{\bm \Pi}_{\bf P,U}$
in Remark $\theREMviS.\theREMviT$(c) must be independent
of the choice of the vectors in ${\bf U}$.
Using the argument of the proof of Proposition $\theLEMivS.\theLEMivT$,
this means that the expression
$$\begin{aligned}
&\Psi_0(\lambda_1) + \Psi_1(\lambda_1)(\lambda_1 - 1) 
+ \Psi_1(\lambda_2)\lambda_1 + \Psi(\lambda_3)\lambda_1 -\cr 
&\hskip7mm-\! \alpha_3 \Big\{ 
\chi_0(\lambda_1,\lambda_2,\lambda_3) \!-\! 
\chi_1(\lambda_1,\lambda_2,\lambda_3)\Big\} \!-\! 
(-1-\alpha_3)  
\chi_1(\lambda_2,\lambda_1,\lambda_3) - \cr
&\hskip7mm-\! \alpha_2 \Big\{ 
\chi_0(\lambda_1,\lambda_3,\lambda_2) \!-\! 
\chi_1(\lambda_1,\lambda_3,\lambda_2)\Big\} \!-\!
(-1-\alpha_2)
\chi_1(\lambda_3,\lambda_1,\lambda_2)
\end{aligned}$$
must be independent of the scalars $\alpha_1,\alpha_2,\alpha_3\in\RR$.
In particular it follows
$\chi_0(\lambda_1,\lambda_2,\lambda_3) \!-\! 
\chi_1(\lambda_1,\lambda_2,\lambda_3) \!-\!
\chi_1(\lambda_3,\lambda_1,\lambda_3) =0$
which completes the proof.

\bigskip
\ADT
\newcounter{REMviiS} \setcounter{REMviiS}{\value{SEC}}
\newcounter{REMviiT} \setcounter{REMviiT}{\value{THM}}
\noindent
{\bf \NUM. Remark.}  (a) According to  
$(\theAFS.\theAFT)$,
affinity invariant tuples are of the form 
\ADT
\newcounter{PSIsS} \setcounter{PSIsS}{\value{SEC}}
\newcounter{PSIsT} \setcounter{PSIsT}{\value{THM}}
$${\bm \Pi}=\big[\Psi,\Psi,\chi_0,\chi_1)\big], \qquad
\chi_0 \big\vert \Delta_3 = 2 \chi_1^{[s]} \big\vert\Delta_3
\leqno(\NUM)$$
with a d-symmetric shape function $\Psi$ and 
the symmetrization $(\theCHInsS.\theCHInsT)$  
of $\chi_1$.

\smallskip
(b) \ From the final argument in the proof of 
Lemma $\thePROPiS.\thePROPiT$ it readily follows
the below converse statement.

\bigskip
\ADT
\newcounter{LEMviiiS} \setcounter{LEMviiiS}{\value{SEC}}
\newcounter{LEMviiiT} \setcounter{LEMviiiT}{\value{THM}}
\noindent
{\bf \NUM. Corollary.} {\it
Assume ${\bm \Pi}=[\Psi,\Psi,\chi_0,\chi_1]$ is an RSD-tuple 
with $(\thePSIsS.\thePSIsT)$.
Then, independently of the choice of the vectors 
${\bf u}_1,{\bf u}_2,{\bf u}_3$ we get
\vskip-7mm
\ADT
\newcounter{LAMS} \setcounter{LAMS}{\value{SEC}}
\newcounter{LAMT} \setcounter{LAMT}{\value{THM}} 
$$\hbox{\Goth F}^{\bm \Pi}_{\bf P \!,U} \lambda_1 =
\sum_{i=1}^3 \Psi(\lambda_i) \lambda_1 +
\chi_1(\lambda_2,\lambda_1,\lambda_3)+\chi_1(\lambda_3,\lambda_1,\lambda_2) .
\leqno(\NUM)$$ 
} 
\DEL{
{\bf Proof.}
Recall that, by writing $f_i = f({\bf p}_i)$,
$g_{ij} = f^\prime({\bf p}_i)[{\bf p}_j-{\bf p}_i]$ 
and 
$\hbox{\Goth F} =\hbox{\Goth F}^{\bm \Pi}_{\bf P \!,U}$
for short,
the function $F=\hbox{\Goth F} f$ 
can be expressed in the form
$$\sum_{i=1}^3 \! \Psi(\lambda_i) f_i + \!\!
\sum_{i,j=1}^3 \! g_{ij} \Psi(\lambda_i) \lambda_j 
-\hskip-3mm \sum_{\hskip-5mm(\ell,m,n)\in{\rm S}_3} \hskip-2mm
{G_\ell {\bf u}_n \over G_n {\bf u}_n} 
\Big\{ \! f_\ell \chi_0(\lambda_\ell,\!\lambda_m,\!\lambda_n) +
g_{\ell m} \chi_1(\lambda_\ell,\!\lambda_m,\!\lambda_n) \!\Big\} .
$$
In particular, if $f=\lambda_1$ we have
$f_1=1$, $f_2=f_3=0$ and $g_{ij} = G_1({\bf p}_j-{\bf p}_i) = 
\lambda_1({\bf p}_j)-\lambda_1({\bf p}_j) =
\delta_{1j}-\delta_{1i}$ $(i,j=1,2,3)$.
Hence $\big($with ordering ${\rm S}_3$ as
$(1,\! 2,\! 3),(2,\! 1,\! 3),(1,\! 3,\! 2),
(3,\! 1,\! 2),(2,\! 3,\! 1),(3,\! 2,\! 1) \big)$,
$$\begin{aligned}
\hbox{\Goth F} &\lambda_1 = \Psi(\lambda_1) - \Psi(\lambda_1)\lambda_2 
- \Psi(\lambda_1)\lambda_3 + \Psi(\lambda_2)\lambda_1 +
\Psi(\lambda_3)\lambda_1 -\cr 
&-\! {G_1 {\bf u}_3 \over G_3 {\bf u}_3} \Big\{ 
\chi_0(\lambda_1,\lambda_2,\lambda_3) \!-\! 
\chi_1(\lambda_1,\lambda_2,\lambda_3)\Big\} \!-\! 
{G_2 {\bf u}_3 \over G_3 {\bf u}_3}  
\chi_1(\lambda_2,\lambda_1,\lambda_3) - \cr
&-\! {G_1 {\bf u}_2 \over G_2 {\bf u}_2} \Big\{ 
\chi_0(\lambda_1,\lambda_3,\lambda_2) \!-\! 
\chi_1(\lambda_1,\lambda_3,\lambda_2) \!-\!
{G_3 {\bf u}_2 \over G_2 {\bf u}_2}
\chi_1(\lambda_3,\lambda_1,\lambda_2) \!-\! 0 \!-\! 0 .
\end{aligned}$$
In view of the relations  
$\lambda_1+\lambda_2+\lambda_3=1$ and 
$(\thePSIsS.\thePSIsT)$,  
we get
\vskip-7mm
$$\begin{aligned}
\hbox{\Goth F} \lambda_1 = \sum_{i=1}^3 \Psi(\lambda_i)\lambda_1 &- 
{G_1 {\bf u}_3 \over G_3 {\bf u}_3} 
\chi_1(\lambda_1,\lambda_2,\lambda_3) - 
{G_2 {\bf u}_3 \over G_3 {\bf u}_3}  
\chi_1(\lambda_2,\lambda_1,\lambda_3) - \cr
&- {G_1 {\bf u}_2 \over G_2 {\bf u}_2}  
\chi_1(\lambda_2,\lambda_1,\lambda_2)\Big\} -
{G_3 {\bf u}_2 \over G_2 {\bf u}_2}
\chi_1(\lambda_3,\lambda_1,\lambda_2) .
\end{aligned}$$
We complete the proof with an application of 
$(\theGUS.\theGUT)$.
}

\DEL{
\ADT
\newcounter{REMixS} \setcounter{REMixS}{\value{SEC}}
\newcounter{REMixT} \setcounter{REMixT}{\value{THM}}
\noindent
{\bf \NUM. Remark.}
(a) \ Given any permutation $(i,j,k)\in{\rm S}_3$, 
the same computation when writing $i_1\!=\! j,i_2\!=\! i,i_3\!=\! k$ in place
of the indices $1,2,3$ yields that
\ADT
\newcounter{LAMiS} \setcounter{LAMiS}{\value{SEC}}
\newcounter{LAMiT} \setcounter{LAMiT}{\value{THM}} 
$$\hbox{\Goth F}^{\bm \Pi}_{\bf P,U} \lambda_j =
\sum_{n=1}^3 \Psi(\lambda_n) \lambda_j +
\chi_1(\lambda_i,\lambda_j,\lambda_k)+\chi_1(\lambda_k,\lambda_j,\lambda_i) .
\leqno(\NUM)$$ 
}

\ADT
\newcounter{PROPviiS} \setcounter{PROPviiS}{\value{SEC}}
\newcounter{PROPviiT} \setcounter{PROPviiT}{\value{THM}}
\noindent
{\bf \NUM. Proposition.} {\it An RSD tuple 
${\bm \Pi}$   
is affinity invariant
if and only if  
it 
is of the form
${\bf \Pi} =[\Psi,\Psi,\chi_0,\chi_1]$  with range shift property, 
$(\thePSIsS.\thePSIsT)$ and
\ADT
\newcounter{CHIstS} \setcounter{CHIstS}{\value{SEC}}
\newcounter{CHIstT} \setcounter{CHIstT}{\value{THM}}
$$\chi_1(t_2,t_1,t_3) \!+\!\chi_1(t_3,t_1,t_2) \!=
t_1 \hskip-7mm \sum_{\hskip5mm(\ell,m,n)\in S_3} \hskip-7mm
\chi_1(t_\ell,t_m,t_n) \quad
\big( (t_1,t_2,t_3)\!\in\!\Delta_3\big) . 
\leqno(\NUM)$$  
} 

\vskip-5mm
{\bf Proof.} Necessity: Assume the affinity invariance of ${\bf \Pi}$. 
According to $\theREMviiS.\theREMviiT$ and 
Corollary $\theLEMviiiS.\theLEMviiiT$,
we can write ${\bm \Pi}=[\Psi,\Psi,\chi_0,\chi_1]$ with
\vskip-7mm
\ADT
\newcounter{LiS} \setcounter{LiS}{\value{SEC}}
\newcounter{LiT} \setcounter{LiT}{\value{THM}}  
$$\lambda_1 = \sum_{i=1}^3 \Psi(\lambda_i)\lambda_1 + 
\chi_1(\lambda_2,\lambda_1,\lambda_3)+\chi_1(\lambda_3,\lambda_1,\lambda_2) .
\leqno(\NUM)$$ 

\vskip-2mm
\noindent
On the other hand 
(Remark $\theREMviS.\theREMviT$(a), Lemma $\thePROPiS.\thePROPiT)$ 
${\bm \Pi}$ has range shift property 
with $(\thePSIsS.\thePSIsT)$ and 
$1 = \sum\limits_{i=1}^3 \Psi(\lambda_i) \!+\!
\chi_0(\lambda_1,\lambda_2,\lambda_3) \!+\! 
\chi_0(\lambda_2,\lambda_3,\lambda_1) \!+\!
\chi_0(\lambda_3,\lambda_1,\lambda_2)$
by $(\theSYMtS.\theSYMtT)$.  Hence, in view of 
$(\thePSIsS.\thePSIsT)$ we conclude that
\vskip-7mm
\ADT
\newcounter{LiiS} \setcounter{LiiS}{\value{SEC}}
\newcounter{LiiT} \setcounter{LiiT}{\value{THM}}  
$$1 = \sum_{i=1}^3 \Psi(\lambda_i)+ 
\hskip-3mm\sum_{(\ell,m,n)\in S_3} \hskip-2mm
\chi_1(\lambda_\ell,\lambda_m,\lambda_n) .
\leqno(\NUM)$$
\vskip-2mm
\noindent
We obtain $(\theCHIstS.\theCHIstT)$ by subtracting
$(\theLiS.\theLiT)$ from $\lambda_1 \cdot (\theLiiS.\theLiiT)$. 

\medskip
Sufficiency: We obtain $(\theLiS.\theLiT)$ meaning the affinity invariance
of ${\bm \Pi}$
by adding the equations 
$0 = \chi_1(\lambda_2,\lambda_1,\lambda_3) \!+\!
\chi_1(\lambda_3,\lambda_1,\lambda_2) -\!
\lambda_1 \hskip-7mm\sum\limits_{\hskip5mm(\ell,m,n)\in S_3} 
\hskip-8mm\chi_1(\lambda_\ell,\lambda_m,\lambda_n)$ \break

\vskip-10mm
\noindent
and
$\lambda_1 \cdot (\theLiiS.\theLiiT)$.

\bigskip
\ADT
\newcounter{CORivS} \setcounter{CORivS}{\value{SEC}}
\newcounter{CORivT} \setcounter{CORivT}{\value{THM}}
\noindent
{\bf \NUM. Corollary.} {\it 
A tuple ${\bf \Pi}=[\Psi_0,\Psi_1,\chi_0,\chi_1]$
with range shift property and 
such that
$\chi_1=\chi_0 \varphi$ for some 
$\varphi \in {\cal C}^1\big({\rm dom}(\chi_0)\big)$ 
is affinity invariant if 
for every $(t_1,t_2,t_3)\in\Delta_3$ we have
\ADT
\newcounter{CHIPHIS} \setcounter{CHIPHIS}{\value{SEC}}
\newcounter{CHIPHIT} \setcounter{CHIPHIT}{\value{THM}}
$$\begin{aligned}
&1 = \varphi(t_1,t_2,t_3) + \varphi(t_2,t_1,t_3) , \cr 
&\chi_0(t_3,t_1,t_2) \varphi(t_3,t_1,t_2) +
\chi_0(t_2,t_1,t_3) \varphi(t_2,t_1,t_3) = \cr
&\hskip7mm = t_1 \Big[  
\chi_0(t_1,t_2,t_3) + \chi_0(t_2,t_3,t_1) + \chi_0(t_2,t_3,t_1)
\Big] .\NUM
\end{aligned}$$
}

\ADT
\newcounter{EXxiiS} \setcounter{EXxiiS}{\value{SEC}}
\newcounter{EXxiiT} \setcounter{EXxiiT}{\value{THM}}
\noindent
{\bf \NUM. Example.} The tuple
${\bm \Pi} \big[ \Phi,\Phi,\chi_0,\chi_0\varphi\big]$ 
is affinity invariant 
with
\vskip-7mm
$$\Phi(t) = t^3 (10 - 15 t + 6 t^2), \ \ 
\chi_0(t_1,t_2,t_3) = 30 t_1^2 t_2^2 t_3, \ \ 
\varphi(t_1,t_2,t_3) = t_2 + {1\over 2} t_3 .$$
Proof. \ The range shift property of ${\bm \Pi}$ 
is established in Example $\theEXivS.\theEXivT$. 
The identities $(\theCHIPHIS.\theCHIPHIT)$
follow with straightforward calculation.

\vskip10mm

\centerline{\bf References}

\medskip

\begin{itemize}

\item[]{[\MDPI]} 
Stach\'o, L.L.;
A Simple Affine-Invariant Spline Interpolation over Triangular Meshes
{\em Mathematics} {\bf 2022}, {\em 10(5)}, 776; 

\smallskip

\item[]{[\STA]} Stach\'o, L.L.;
Locally Generated Polynomial C1-Splines over Triangular Meshes.
{\em Miskolc Mathematical Notes} {\bf 2022}, {\em 23(2)}, 897--911.

\smallskip

\item[][\Zlam]
Zl\'amal, M.; Z\v eni\v sek, 
A. Mathematical aspects of the FEM.
In {\em Technical Physical and Mathematical Priciples of the FEM}; Kola\v z, V., Kratoci-Ivil, F., Zlamal, M., Zenisek, A., 
Eds.; Publising House: Praha, Czech Republic, {\bf 1971}; pp. 15--39.

\smallskip

\item[][\Serg]
I.V. Sergienko, O.M. Lytvyn, O.O. Lytvyn and O.I. Denisova;
Explicit formulas for interpolating splines of degree 5 on the triangle,
{\em Cybernetics and Systems Analysis} {\bf 2014}, {\em 50(5)}, 670--678.

\smallskip

\item[][\Hahm]
Hahmann, S.; Bonneau, G.-P. 
Triangular $G^1$ interpolation by 4-splitting domain triangles.
{\em Comput. Aided Geom. Des.} {\bf 2000}, {\em 17}, 731--757.

\smallskip

\item[][\Cao]
Cao, J.; Chen, Z.; Wei, X.; Zhang, Y.J. 
A finite element framework based on bivariate simplex splines
on triagle configurations.
{\em Comput. Methods Appl. Mech. Eng.} {\bf 2019}, {\em 357}, 112598.

\smallskip

\item[][\WHIT]
Whitney, Hassler; 
Analytic extensions of differentiable functions, 
{\em Transactions of the Amer. Math. Soc.} {\bf 1934}, 
{\em 36(1)}, 63--89. 

\smallskip

\item[][\Berger]
Berger, Marcel; {Geometry I}; Springer, New York, {\bf 2009}.

\end{itemize}

\end{document}